\documentclass{revtex4}
\usepackage{graphicx,amsmath,amssymb,mathrsfs}
\usepackage{epsfig}
\usepackage{caption}
\usepackage{mathtools}
\usepackage{color}
\usepackage{subfig}
\begin{document}

\title{Stochastic Differential Equations Driven by Deterministic Chaotic Maps: Analytic Solutions of the Perron-Frobenius Equation}

\author{Griffin Williams and Christian Beck}

\affiliation{School of Mathematical Sciences, Queen Mary University of London, Mile End Road, London E1 4NS, UK}

\begin{abstract}
We consider discrete-time dynamical systems with a linear relaxation dynamics that are driven by deterministic chaotic forces.  By perturbative expansion in a small time scale parameter, we derive from the Perron-Frobenius equation the corrections to ordinary Fokker-Planck equations in leading order of the time scale separation parameter.  We present analytic solutions to the equations for the example of driving forces generated by $N$-th order Chebychev maps.  The leading order corrections are universal for $N\geq 4$ but different for $N=2$ and $N=3$. We also study
diffusively coupled Chebychev maps as driving forces, where strong correlations may prevent convergence to
Gaussian limit behavior.
\end{abstract}

\maketitle

\section{Introduction}
A fundamental problem lying at the roots of the foundations of statistical mechanics is the question of how Brownian motion and Langevin processes can arise from deterministic dynamical systems with strong mixing properties.
 This problem has been studied by many authors \cite{Bill}-\cite{new5}. Due to the nature of the mixing property
              it is clear that a large phase-space dimension is not necessary to obtain `random behavior' of Brownian motion or Ornstein-Uhlenbeck type \cite{Bill}, rather it is sufficient to consider low-dimensional mixing deterministic maps (with random initial conditions) and perform a suitable scaling limit, meaning looking at the ensemble of trajectories in a suitable rescaled way \cite{Roep}-\cite{Beck5}. This is mathematically described by a time scale parameter $\tau$
which approaches 0. This parameter describes the ratio of the fast time scale underlying the deterministic chaotic
driving force and the slow time scale corresponding to the relaxation time of the system.

In the simplest case, namely the case of a linear Langevin equation
driven by chaotic `noise', the relevant class of dynamical systems
where these concepts can be illustrated nicely, and mathematically treated rigorously,
is given by maps of Kaplan-Yorke type \cite{Kaplan}, sometimes also called {\em maps of linear Langevin type} \cite{Beck6} or {\em skew products} \cite{melbourne}.  They can be regarded as the deterministic chaotic analogue of a linear Langevin equation \cite{Kampen} where the Gaussian white noise is replaced by a chaotic dynamics.
In general, maps of linear Langevin type generate complicated (non-Gaussian, non-Markovian) stochastic processes,
which, however, reduce to ordinary Gaussian behavior of Ornstein-Uhlenbeck type in the scaling limit $\tau \to 0$ mentioned above, provided the driving map
has sufficiently strong mixing behavior.

Of course, for small but finite $\tau$ there are corrections to the Gaussian limit behavior due to the underlying deterministic dynamics, which always has nontrivial correlations on a small scale, as it
  is deterministic on a small scale. In this paper we deal with the general question of what these corrections
  to Gaussian Langevin behavior are for general mixing driving forces,
for small but finite $\tau$.  The driving dynamics is assumed to be generated by a map $T$ that has the so-called $\varphi$-mixing property, a property that is sufficient to guarantee convergence to
a Brownian motion process when the iterates are summed up and properly rescaled.
We start from the Perron-Frobenius equation and then derive, by a perturbative
expansion in $\sqrt{\tau}$,
         a set of equations that is obeyed by the corrections to the Gaussian limit behavior. In 0-th order
         we get from this the Fokker-Planck equation, but our main goal
         is to solve the next-to-leading order equations. These turn out to be complicated functional equations.
         In this paper we show that we can explicitly solve the equations if the chaotic driving force is a
         Chebychev map of $N$-th order, with $N$ arbitrary.

More complicated cases arise if the linear relaxation dynamics is driven by coupled maps
\cite{kaneko,buni}. A dynamics
driven by diffusively coupled Chebychev maps has applications in quantum field theory \cite{Beck7}-\cite{SM},
providing a discrete dynamical version of the Parisi-Wu approach of stochastic quantization.
The Perron-Frobenius equation cannot be solved analytically in this more complicated case, but still
we will present some numerical results indicating how typical densities look like in this case,
in particular how the kurtosis depends on the coupling strength of the coupled map
lattice for small values of $\tau$.

This paper is organized as follows:

In section 2 we perturbatively expand the Perron-Frobenius equation
for the general type of dynamical systems considered in this paper. In section 3 we integrate out the
chaotic forcing degrees
of freedom and proceed to marginal densities (which correspond to the velocity of the kicked particle
in the physical picture). In the following sections we then present the solution of the equations for general
$N$-th order Chebychev maps $T_{N}$ as discrete-time driving forces. Each case $N=2,3,4$ requires a considerable amount of
calculations. We start with the case $N=2$ in section 4, and then solve the
functional equations for $N=3$ in section 5, and $N\geq 4$ in section 6.  Technical
details of the calculations can be found in the appendix.

In section 7 we will then finally compare our analytical predictions with numerical results.
For the different $N$ we compare our analytically obtained results with numerically obtained histograms of
iterates for different values of the time scale parameter $\tau$.
Excellent agreement between analytical and numerical results is obtained. In this section we will also present numerical results for cases that are
not analytically tractable, such as the coupled Chebychev maps as driving forces as mentioned above.
Our concluding remarks are given in section 8.

\section{Perturbative Expansion of the Perron-Frobenius Equation in the time scale parameter $\sqrt{\tau}$}
We start from a dynamical system of linear Langevin type \cite{Roep,Beck1,Shim1,Beck5,new2,new4,new5,Yalcin,Schlogl,HC1,HC2}
\begin{equation}
f:
        \begin{array}{ll}
                x_{n+1} = T \left(x_{n}\right)\\
                y_{n+1} = \lambda y_{n} + \tau^{1/2}x_{n}
        \end{array}
\end{equation}
Here $\lambda\in(0,1)$ and $\tau >0$ are parameters, and $T: X\rightarrow X$ is some mapping that has strong mixing properties such that convergence to a Gaussian stochastic process $Y(t)$ is guaranteed in the limit $\tau\rightarrow0, \ \lambda\rightarrow1, \ t=n\tau$ finite.  For example, the so-called $\varphi$-mixing property is sufficient.  We recall here that a stationary sequence $\lbrace\xi_{n}\rbrace$ of random variables is said to be $\varphi$-mixing if the function $\varphi(n)$, defined by
\begin{equation} \nonumber
	\varphi(n)= \operatorname*{sup}_{k\geqslant0} \ \operatorname*{\operatorname*{sup}_{A\in \mathscr{B}_{k}}}_{B\in \mathscr{B}^{k}} \ \vert m(T^{-n}A\vert B)-m(A)\vert \ \ (n\geq0)
\end{equation}
has the property $\varphi(n)\rightarrow0 \ (n\rightarrow\infty)$. Here $m$ denotes a suitable
 (natural) invariant measure of $T$. Let $\mathscr{B}^{n}$ denote the $\sigma$-algebra generated by the random variables $\xi_{0},\dots,\xi_{n-1}$, and $\mathscr{B}_{n}$ denote the $\sigma$-algebra generated by $\xi_{n},\xi_{n+1}\dots$.  A map $T$ is said to be $\varphi$-mixing if it admits a generating partition $\xi$ such that the sequence $\xi_{n}=\xi\circ T^{n}$ is $\varphi$-mixing \cite{Bill,Roep}.  As a simple
 example, any map conjugated to a Bernoulli shift is $\varphi$-mixing (with a vanishing $\varphi$-function).

 The map $f$ is a skew product extension of $T$. Physically it is obtained by integration from the following deterministic chaotic analogue of a Langevin equation,
which formally describes the velocity $Y$ of a damped particle under deterministic chaotic kicks $x_n$:
\begin{equation}
	\dot{Y}=-\gamma Y + \tau^{1/2}\sum_{n=1}^{\infty} x_{n-1}\delta(t-n\tau)
\end{equation}
\begin{equation}
	x_{n+1}=T(x_{n})
\end{equation}
The map $T$ determines the time evolution of the kicks. $Y(t)$ can be regarded as the velocity of a kicked damped particle.  Equation $(1)$ describes the stroboscopic time evolution of $y_{n}:=Y(n\tau +0)$.  The damping constant $\gamma >0$ and the time difference $\tau$ between kicks are related to the parameter $\lambda$ by $\lambda = e^{-\gamma\tau}$.  For concreteness we will deal with one-dimensional driving forces, although a similar analysis can be performed for higher-dimensional cases.  The inverse map is given by
\begin{equation}
f^{-1}:
	\begin{array}{ll}
		x_{n}=T^{-1}(x_{n+1})\\
		y_{n}=\lambda^{-1}(y_{n+1}-\tau^{1/2}T^{-1}(x_{n+1}))
	\end{array}
\end{equation}
In general, there are several preimages $T^{-1}$.  The absolute value of the Jacobi determinant of $f$ is
\begin{equation}
	\vert det Df \vert = \lambda\vert T'(x)\vert
\end{equation}
The Perron-Frobenius equation
\begin{equation}
	\rho_{n+1}(x',y') = \sum_{(x,y)\in f^{-1}(x',y')} \frac{\rho_{n}(x,y)}{\lambda\vert T'(x)\vert}
\end{equation}
governs the time evolution of probability densities $\rho_{n}(x,y)$; it can be written as a sum over the preimages of $T$:
\begin{equation}
	\lambda\rho_{n+1}(x',y)=\sum_{x\in T^{-1}(x')}\frac{1}{\vert T'(x)\vert}\rho_{n}(x,\lambda^{-1}(y-\tau^{1/2}x))
\end{equation}
Equation $(7)$, with $\lambda = e^{-\gamma\tau}$, is the starting point of our consideration.  We will expand it with respect to the parameter $\tau^{1/2}$, which is supposed to be small. Our main assumption is that such an expansion does make sense for suitable classes
of maps $T$ and leads to a convergent series.  We cannot prove this in general, but we will consider concrete examples of strongly mixing maps ($N$-th order Chebychev maps) in this paper where indeed we show that the equations can be solved term by term in $\tau^{1/2}$, leading to a series expansion that correctly describes the scaling behavior of the invariant density for small values of $\tau$.  In our perturbative analysis we will include all terms up to fourth order in $\tau^{1/2}$.  To simplify the notation, we will not explicitly write down any term of $O(\tau^{5/2})$ or higher, but suppress it in all equations.  Our numerics indeed indicates that higher order terms are negligible if $\tau$ is small enough (details of the numerics in Section VII). Without restriction of generality we set $\gamma=1$, obtaining up to fourth order in $\tau^{1/2}$
\begin{equation}
	\lambda = e^{-\tau} = 1-\tau +\frac{1}{2}\tau^{2}
\end{equation}
Since
\begin{equation}
	\lambda^{-1}(y-\tau^{1/2}x)=(1+\tau +\frac{1}{2}\tau^{2})(y-\tau^{1/2}x)
\end{equation}
\begin{equation}
	=y-\tau^{1/2}x+\tau y-\tau^{3/2}x+\frac{1}{2}\tau^{2}y
\end{equation}
we obtain by Taylor expansion
\begin{equation}
\begin{aligned}
	\rho_{n}(x,\lambda^{-1}(y-\tau^{1/2}x)) & \\
	& ={} \rho_{n}(x,y)+\left(-\tau^{1/2}x+\tau y-\tau^{3/2}x+\frac{1}{2}\tau^{2}y\right)\frac{\partial}{\partial y}\rho_{n}(x,y) \\
	& +\frac{1}{2}(\tau x^{2}+\tau^{2}y^{2}-2\tau^{3/2}xy+2\tau^{2}x^{2})\frac{\partial^{2}}{\partial y^{2}}\rho_{n}(x,y) \\
	& +\frac{1}{6}(-\tau^{3/2}x^{3}+3\tau^{2}x^{2}y)\frac{\partial^{3}}{\partial y^{3}}\rho_{n}(x,y)+\frac{1}{24}\tau^{2}x^{4}\frac{\partial^{4}}{\partial y^{4}}\rho_{n}(x,y) \\
	& =\rho_{n}(x,y)+\tau^{1/2}\left[-x\frac{\partial}{\partial y}\rho_{n}(x,y)\right] \\
	& +\tau\left[y\frac{\partial}{\partial y}\rho_{n}(x,y)+\frac{1}{2}x^{2}\frac{\partial^{2}}{\partial y^{2}}\rho_{n}(x,y)\right] \\
	& +\tau^{3/2}\left[-x\frac{\partial}{\partial y}\rho_{n}(x,y)-xy\frac{\partial^{2}}{\partial y^{2}}\rho_{n}(x,y)-\frac{1}{6}x^{3}\frac{\partial^{3}}{\partial y^{3}}\rho_{n}(x,y)\right] \\
	& +\tau^{2}\biggl[\frac{1}{2}y\frac{\partial}{\partial y}\rho_{n}(x,y)+\frac{1}{2}y^{2}\frac{\partial^{2}}{\partial y^{2}}\rho_{n}(x,y)+x^{2}\frac{\partial^{2}}{\partial y^{2}}\rho_{n}(x,y) \\
	& +\frac{1}{2}x^{2}y\frac{\partial^{3}}{\partial y^{3}}\rho_{n}(x,y)+\frac{1}{24}x^{4}\frac{\partial^{4}}{\partial y^{4}}\rho_{n}(x,y)\biggr]
\end{aligned}
\end{equation}
Let us introduce a continuous-time suspension $\rho(x,y,t)$ defined by
\begin{equation}
	\rho_{n}(x,y)=\rho(x,y,t) \ \ \ \ (t=n\tau)
\end{equation}
By such a continuous-time suspension we mean a function $\rho(x,y,t)$ depending on a continuous time variable $t$ such that at stroboscopic times $t=n\tau$ it coincides with $\rho_{n}(x,y)$.  Since we keep $\tau$ small but finite, there are infinitely many such smooth functions $\rho(x,y,t)$.  At the present stage, we need not fix the function $\rho(x,y,t)$ for time values $t$ other than
discrete values $n\tau$.  The following considerations are valid for \textit{any} smooth suspension $\rho(x,y,t)$ satisfying Eq. $(12)$.  A Taylor expansion yields
\begin{equation} \nonumber
\begin{aligned}
	\rho_{n+1}(x,y)={} & \rho(x,y,n\tau +\tau) \\
	 & =\rho(x,y,n\tau)+\tau\frac{\partial}{\partial t}\rho(x,y,n\tau)+\frac{1}{2}\tau^{2}\frac{\partial^{2}}{\partial t^{2}}\rho(x,y,n\tau)
\end{aligned}
\end{equation}
Hence
\begin{equation}
\begin{aligned}
	\lambda\rho_{n+1}(x',y) & \\
	& ={} \left(1-\tau +\frac{1}{2}\tau^{2}\right)\left[\rho(x',y,t)+\tau\frac{\partial}{\partial t}\rho(x',y,t)+\frac{1}{2}\tau^{2}\frac{\partial^{2}}{\partial t^{2}}\rho(x',y,t)\right] \\
	& =\rho(x',y,t)+\tau\left[-\rho(x',y,t)+\frac{\partial}{\partial t}\rho(x',y,t)\right] \\
	& +\tau^{2}\left[\frac{1}{2}\rho(x',y,t)-\frac{\partial}{\partial t}\rho(x',y,t)+\frac{1}{2}\frac{\partial^{2}}{\partial t^{2}}\rho(x',y,t)\right] \ \ \ \ (t=n\tau)
\end{aligned}
\end{equation}
The function $\rho(x,y,t)$ still depends on $\lambda$ and thus on $\tau^{1/2}$, since for stroboscopic times $t=n\tau$ it is the density $\rho(x,y)$ of the map $(1)$ which explicitly depends on $\lambda$.  In the spirit of van Kampen's $\Omega$-expansion \cite{Kampen}, we now make the following ansatz:
\begin{equation}
\begin{aligned}
	\rho(x,y,t)={} & \varphi(x,y,t)+\tau^{1/2}a(x,y,t)+\tau b(x,y,t) \\
	& +\tau^{3/2}c(x,y,t)+\tau^{2}d(x,y,t),
\end{aligned}
\end{equation}
where the functions $\varphi, a, b, c,$ and $d$ are assumed to be independent of $\tau$.
 In our expansion all functions $\varphi, a,b,c,d$ are assumed to be smooth such that they can be differentiated as often as required.
 For the concrete examples of maps $T$ studied in this paper, we will indeed find such smooth differentiable solutions to the equations in the later sections.

     By putting the ansatz $(14)$ into Eqs. $(11)$ and $(13)$ and comparing different powers of $\tau^{1/2}$, one finally obtains the following five coupled functional equations for $\varphi, a, b, c,$ and $d$:
\begin{equation}
             \varphi(x',y,t)= \sum_{x\in T^{-1}(x')}\frac{1}{\vert T'(x)\vert}\varphi(x,y,t)
\end{equation}
\begin{equation}
			 a(x',y,t)= \sum_{x\in T^{-1}(x')}\frac{1}{\vert T'(x)\vert}\biggl[a(x,y,t)-x\frac{\partial}{\partial y}\varphi(x,y,t)\biggr]
\end{equation}
\begin{equation}
\begin{aligned}
b(x',y,t) ={} & \sum_{x\in T^{-1}(x')}\frac{1}{\vert T'(x)\vert}\biggl[b(x,y,t)-x\frac{\partial}{\partial y}a(x,y,t) \\
			& +\frac{\partial}{\partial y}(y\varphi(x,y,t))+\frac{1}{2}x^{2}\frac{\partial^{2}}{\partial y^{2}}\varphi(x,y,t)-\frac{\partial}{\partial t}\varphi(x,y,t)\biggr]
\end{aligned}
\end{equation}
\begin{equation}
\begin{aligned}
c(x',y,t) ={} & \sum_{x\in T^{-1}(x')}\frac{1}{\vert T'(x)\vert}\biggl[c(x,y,t)-x\frac{\partial}{\partial y}b(x,y,t)\\
		& +\frac{\partial}{\partial y}(ya(x,y,t))+\frac{1}{2}x^{2}\frac{\partial^{2}}{\partial y^{2}}a(x,y,t)-\frac{\partial}{\partial t}a(x,y,t) \\
		& -2x\frac{\partial}{\partial y}\varphi(x,y,t)+x\frac{\partial}{\partial y}\frac{\partial}{\partial t}\varphi(x,y,t)-xy\frac{\partial^{2}}{\partial y^{2}}\varphi(x,y,t) \\
		& -\frac{1}{6}x^{3}\frac{\partial^{3}}{\partial y^{3}}\varphi(x,y,t)\biggr]
\end{aligned}
\end{equation}
\begin{equation}
\begin{aligned}
d(x',y,t) ={} & \sum_{x\in T^{-1}(x')}\frac{1}{\vert T'(x)\vert}\bigg\{d(x,y,t)-x\frac{\partial}{\partial y}c(x,y,t)+\frac{\partial}{\partial y}(yb(x,y,t)) \\
			& +\frac{1}{2}x^{2}\frac{\partial^{2}}{\partial y^{2}}b(x,y,t)-\frac{\partial}{\partial t}b(x,y,t) \\
			& +x\left[-2 +\frac{\partial}{\partial t}-y\frac{\partial}{\partial y}-\frac{1}{6}x^{2}\frac{\partial^{2}}{\partial y^{2}}\right]\frac{\partial}{\partial y}a(x,y,t) \\
			& +\frac{\partial}{\partial y}(y\varphi(x,y,t))+\frac{1}{2}x^{2}\frac{\partial^{2}}{\partial y^{2}}\varphi(x,y,t)-\frac{\partial}{\partial t}\varphi(x,y,t) \\
			& -\frac{\partial}{\partial t}\left[\frac{\partial}{\partial y}(y\varphi(x,y,t))+\frac{1}{2}x^{2}\frac{\partial^{2}}{\partial y^{2}}\varphi(x,y,t)-\frac{\partial}{\partial t}\varphi(x,y,t)\right] \\
			& +\frac{1}{2}y\frac{\partial}{\partial y}\varphi(x,y,t)+\frac{1}{2}y^{2}\frac{\partial^{2}}{\partial y^{2}}\varphi(x,y,t) \\
			& +x^{2}\left[1+\frac{1}{2}y\frac{\partial}{\partial y}+\frac{1}{24}x^{2}\frac{\partial^{2}}{\partial y^{2}}\right]\frac{\partial^{2}}{\partial y^{2}}\varphi(x,y,t) \\
			& -\frac{1}{2}\varphi(x,y,t)+\frac{\partial}{\partial t}\varphi(x,y,t)-\frac{1}{2}\frac{\partial^{2}}{\partial t^{2}}\varphi(x,y,t)\bigg\}
\end{aligned}
\end{equation}
Notice that the last term on the right-hand side of Eq. $(17)$ is a Fokker-Planck operator with \textit{non-constant} variance $x^{2}$, similar in spirit as for the superstatistics approach \cite{super1,super2}.  Note that in the superstatistics approach one studies spatio-temporal systems with inhomogeneous temperature distributions $T$, meaning the diffusion constant $D\sim kT$ is varying in space and/or time.  The term $\frac{1}{2}x^{2}$ in front of the term $\frac{\partial^{2}}{\partial y^{2}} \varphi(x,y,t)$ in Eq. (17) can formally be interpreted as a varying diffusion constant of a local Fokker-Planck equation for the function $\varphi(x,y,t)$, and later we will integrate over all possible $x$,
similar in spirit to the superstatistics approach, where an ensemble of diffusion constants
is considered \cite{super1,super2,nature-e}.
 Equation $(17)$ is a kind of combination of a Perron-Frobenius equation for the map $T$ with a Fokker-Planck equation, where the variance is not constant, but given by $x^{2}$.
\section{Equations for marginal velocity distributions}
In physical measurements of a particle that moves under deterministic chaotic kicks $x_{n}$, one is interested in the probability distribution of the velocity $y_{n}$ of the particle, and not in the joint probability distribution of kick strength $x_{n}$ and velocity $y_{n}$.  Therefore, from a physical point of view it makes sense to proceed to the marginal distribution
\begin{equation}
	p(y,t)=\int dx \ \rho(x,y,t)
\end{equation}
of the $y$ variable, obtained by integration over all possible $x$ values.   In the physical picture, $p (y,t)$ describes the velocity distribution of the particle, which is kicked by deterministic chaotic iterates $x_{n}$, whereas $\rho(x,y,t)$ has no direct physical interpretation. For $t \to \infty$,
 the stationary distribution is obtained, denoted by $p (y)$. The integration over $x$ indeed yields a simplification in our problem, transforming the coupled functional equations of the previous section into coupled differential equations.

Let us use the notation
\begin{equation}
	p_{0}(y,t) =\int dx \ \varphi(x,y,t)
\end{equation}
\begin{equation}
	\alpha(y,t) =\int dx \ a(x,y,t)
\end{equation}
\begin{equation}
	\beta(y,t) =\int dx \ b(x,y,t)
\end{equation}
\begin{equation}
	\gamma(y,t) =\int dx \ c(x,y,t)
\end{equation}
\begin{equation}
	\delta(y,t) =\int dx \ d(x,y,t)
\end{equation}
From Eqs.(15)-(19) we obtain that the marginal functions satisfy
\begin{equation}
	\frac{\partial}{\partial y}\int dx \ x\varphi(x,y,t)=0
\end{equation}
\begin{equation}
	\frac{\partial}{\partial y}\int dx \ xa(x,y,t)=\frac{\partial}{\partial y}(yp_{0}(y,t))+\frac{1}{2}\frac{\partial^{2}}{\partial y^{2}}\int dx \ x^{2}\varphi(x,y,t)-\frac{\partial}{\partial t}p_{0}(y,t)
\end{equation}
\begin{equation}
\begin{aligned}
	\frac{\partial}{\partial y}\int dx \ xb(x,y,t) ={} & \frac{\partial}{\partial y}(y\alpha(y,t))+\frac{1}{2}\frac{\partial^{2}}{\partial y^{2}}\int dx \ x^{2}a(x,y,t)-\frac{\partial}{\partial t}\alpha(y,t) \\
	& +\int dx \ x\left[-2+\frac{\partial}{\partial t}-y\frac{\partial}{\partial y}-\frac{1}{6}x^{2}\frac{\partial^{2}}{\partial y^{2}}\right]\frac{\partial}{\partial y}\varphi(x,y,t)
\end{aligned}
\end{equation}
\begin{equation}
\begin{aligned}
	\frac{\partial}{\partial y}\int dx \ xc(x,y,t) ={} & \frac{\partial}{\partial y}(y\beta(y,t))+\frac{1}{2}\frac{\partial^{2}}{\partial y^{2}}\int dx \ x^{2}b(x,y,t)-\frac{\partial}{\partial t}\beta(y,t) \\
		& -\int dx \ x\left[2-\frac{\partial}{\partial t}+y\frac{\partial}{\partial y}+\frac{1}{6}x^{2}\frac{\partial^{2}}{\partial y^{2}}\right]\frac{\partial}{\partial y}a(x,y,t) \\
		& +\left(1-\frac{\partial}{\partial t}\right)\biggl[\frac{\partial}{\partial y}(yp_{0}(y,t)) \\
		& +\frac{1}{2}\frac{\partial^{2}}{\partial y^{2}}\int dx \ x^{2}\varphi(x,y,t)-\frac{\partial}{\partial t}p_{0}(y,t)\biggr] \\
		& +\left[\frac{1}{2}y\frac{\partial}{\partial y}+\frac{1}{2}y^{2}\frac{\partial^{2}}{\partial y^{2}}-\frac{1}{2}+\frac{\partial}{\partial t}-\frac{1}{2}\frac{\partial^{2}}{\partial t^{2}}\right]p_{0}(y,t) \\
		& +\left(\frac{\partial^{2}}{\partial y^{2}}+\frac{1}{2}y\frac{\partial^{3}}{\partial y^{3}}\right)\int dx \ x^{2}\varphi(x,y,t) \\
		& +\frac{1}{24}\frac{\partial^{4}}{\partial y^{4}}\int dx \ x^{4}\varphi(x,y,t)
\end{aligned}
\end{equation}
Let us first deal with the zeroth-order term $\varphi(x,y,t)$.  Equation $(15)$ is solved by any function $\varphi(x,y,t)$ of the form
\begin{equation}
	\varphi(x,y,t)=h(x)p_{0}(y,t)
\end{equation}
where $h(x)$ is the natural invariant density of the $T$-dynamics.  To obtain a compact notation, we will use the notation $\langle\dots \rangle$ for expectations with respect to $h(x)$.  Equation $(26)$ then implies
\begin{equation}
	\frac{\partial}{\partial y}p_{0}(y,t)\int dx \ xh(x)=\langle x\rangle\frac{\partial}{\partial y}p_{0}(y,t)=0
\end{equation}
In case that $\rho_{0}(y,t)\neq const_{y}$ this means $\langle x\rangle=0$, which shows that the ansatz $(14)$ makes sense for maps with vanishing average only.  Equation $(27)$ becomes
\begin{equation}
	\frac{\partial}{\partial y}(yp_{0}(y,t))+\frac{1}{2}\langle x^{2}\rangle\frac{\partial^{2}}{\partial y^{2}}p_{0}(y,t)-\frac{\partial}{\partial t}p_{0}(y,t)=\frac{\partial}{\partial y}\int dx \ xa(x,y,t)
\end{equation}
This is a kind of inhomogeneous Fokker-Planck equation with a source term $(\partial /\partial y)\int dx xa(x,y,t)$.  It reduces to a Fokker-Planck equation for the case that $(\partial /\partial y)\int dx xa(x,y,t)=0$.  If $(\partial /\partial y)\int dx x a(x,y,t)$ is proportional to $(\partial^{2}/\partial y^{2})\rho_{0}(y,t)$, we also obtain a Fokker-Planck equation, but with a different diffusion constant.  Hence, in order to determine $\rho_{0}(y,t)$, we have to determine $\int dx xa(x,y,t)$ with the help of Eq. $(16)$.

The equations simplify considerably for particular choices of mappings
     $T$ that have the property that both the mapping $T$ as well as its natural invariant density $h$ are symmetric:
\begin{equation}
	T(x)=T(-x)
\end{equation}
\begin{equation}
	h(x)=h(-x)
\end{equation}
We call these types of maps {\em double symmetric.} Examples are maps conjugated to even Chebyshev polynomials.  For double symmetric maps we have
\begin{equation}
	\sum_{x\in T^{-1}(x')}\frac{1}{\vert T'(x)\vert}xh(x)\frac{\partial}{\partial y}p_{0}(y,t)=0
\end{equation}
since for each $x,-x$ is also a preimage.  Hence Eq. $(16)$ reduces to
\begin{equation}
	a(x',y,t)=\sum_{x\in T^{-1}(x')}\frac{1}{\vert T'(x)\vert}a(x,y,t)
\end{equation}
Notice that this is just the same equation as the one satisfied by $\varphi$.  It is solved by a function of the form
\begin{equation}
	a(x,y,t)=h(x)\alpha(y,t)
\end{equation}
Moreover, we obtain
\begin{equation}
\begin{aligned}
	b(x',y,t) ={} & \sum_{x\in T^{-1}(x')}\frac{1}{\vert T'(x)\vert}(b(x,y,t)+h(x)(\frac{\partial}{\partial y}(yp_{0}(y,t)) \\
		& +\frac{1}{2}x^{2}\frac{\partial^{2}}{\partial y^{2}}p_{0}(y,t)-\frac{\partial}{\partial t}p_{0}(y,t)))
\end{aligned}
\end{equation}
\begin{equation}
\begin{aligned}
	c(x',y,t) ={} & \sum_{x\in T^{-1}(x')}\frac{1}{\vert T'(x)\vert}(c(x,y,t)-x\frac{\partial}{\partial y}b(x,y,t) \\
		& +h(x)\left[\frac{\partial}{\partial y}(y\alpha(y,t))+\frac{1}{2}x^{2}\frac{\partial^{2}}{\partial y^{2}}\alpha(y,t)-\frac{\partial}{\partial t}\alpha(y,t)\right])
\end{aligned}
\end{equation}
\begin{equation}
\begin{aligned}
	d(x',y,t) ={} & \sum_{x\in T^{-1}(x')}\frac{1}{\vert T'(x)\vert}(d(x,y,t)-x\frac{\partial}{\partial y}c(x,y,t)+\frac{\partial}{\partial y}(yb(x,y,t)) \\
		& +\frac{1}{2}x^{2}\frac{\partial^{2}}{\partial y^{2}}b(x,y,t)-\frac{\partial}{\partial t}b(x,y,t) \\
		& +h(x)(\left(1-\frac{\partial}{\partial t}\right)(\frac{\partial}{\partial y}(yp_{0}(y,t)) \\
		& +\frac{1}{2}x^{2}\frac{\partial^{2}}{\partial y^{2}}p_{0}(y,t)-\frac{\partial}{\partial t}p_{0}(y,t)) \\
		& +\frac{1}{2}y\frac{\partial}{\partial y}p_{0}(y,t)+\frac{1}{2}y^{2}\frac{\partial^{2}}{\partial y^{2}}p_{0}(y,t) \\
		& +x^{2}\left(1+\frac{1}{2}y\frac{\partial}{\partial y}+\frac{1}{24}x^{2}\frac{\partial^{2}}{\partial y^{2}}\right)\frac{\partial^{2}}{\partial y^{2}}p_{0}(y,t) \\
		& -\frac{1}{2}p_{0}(y,t)+\frac{\partial}{\partial t}p_{0}(y,t)-\frac{1}{2}\frac{\partial^{2}}{\partial t^{2}}p_{0}(y,t)))
\end{aligned}
\end{equation}
The integrated equations reduce to
\begin{equation}
	0=\frac{\partial}{\partial y}(yp_{0}(y,t))+\frac{1}{2}\langle x^{2}\rangle\frac{\partial^{2}}{\partial y^{2}}p_{0}(y,t)-\frac{\partial}{\partial t}p_{0}(y,t)
\end{equation}
\begin{equation}
	 \frac{\partial}{\partial y}\int dx \ xb(x,y,t)=\frac{\partial}{\partial y}(y\alpha(y,t))+\frac{1}{2}\langle x^{2}\rangle\frac{\partial^{2}}{\partial y^{2}}\alpha(y,t)-\frac{\partial}{\partial t}\alpha(y,t)
\end{equation}
\begin{equation}
\begin{aligned}
	\frac{\partial}{\partial y}\int dx \ xc(x,y,t) ={} & \frac{\partial}{\partial y}(y\beta(y,t))+\frac{1}{2}\frac{\partial^{2}}{\partial y^{2}}\int dx \ x^{2}b(x,y,t)-\frac{\partial}{\partial t}\beta(y,t) \\
		& +\biggl[\frac{1}{2}y\frac{\partial}{\partial y}+\frac{1}{2}y^{2}\frac{\partial^{2}}{\partial y^{2}}+\langle x^{2}\rangle\frac{\partial^{2}}{\partial y^{2}}+\frac{1}{2}\langle x^{2}\rangle y\frac{\partial^{3}}{\partial y^{3}} \\
		& +\frac{1}{24}\langle x^{4}\rangle\frac{\partial^{4}}{\partial y^{4}}-\frac{1}{2}+\frac{\partial}{\partial t}-\frac{1}{2}\frac{\partial^{2}}{\partial t^{2}}\biggr]p_{0}(y,t)
\end{aligned}
\end{equation}
Notice that Eq. $(41)$ is the Fokker-Planck equation.  We obtain the result that for double-symmetric maps the diffusion constant is always given by $\langle x^{2}\rangle$.  Equation $(42)$ is the same type as Eq. $(32)$; it is an inhomogeneous Fokker-Planck equation, but now - due to the condition of double symmetry - it is satisfied by $\alpha(y,t)$ rather that $p_{0}(y,t)$.

Subtracting the zero given by Eq. $(41)$ from the right-hand side of Eq. $(38)$, we may also write Eq. $(38)$ as
\begin{equation}
b(x',y,t)=\sum_{x\in T^{-1}(x')}\frac{1}{\vert T'(x)\vert}\left[b(x,y,t)+h(x)\frac{1}{2}(x^{2}-\langle x^{2}\rangle)\frac{\partial^{2}}{\partial y^{2}}p_{0}(y,t)\right]
\end{equation}
Similarly, Eq. $(39)$ can be written as
\begin{equation}
\begin{aligned}
c(x',y,t) ={} & \sum_{x\in T^{-1}(x')}\frac{1}{\vert T'(x)\vert}\bigg\{c(x,y,t)-x\frac{\partial}{\partial y}b(x,y,t) \\
	& +h(x)\left[\frac{\partial}{\partial y}\int dx \ xb(x,y,t)+\frac{1}{2}(x^{2}-\langle x^{2}\rangle)\frac{\partial^{2}}{\partial y^{2}}\alpha(y,t)\right]\bigg\}
\end{aligned}
\end{equation}
\section{Explict solution for $N=2$, $T_{2}(x)=2x^{2}-1$}
After having written down the general equations satisfied by the higher order corrections to the Fokker-Planck equation for $\varphi$-mixing chaotic maps as driving forces, we now show that the equations can be explicitly solved for Chebyshev maps of $N$-th order ($N \geq 2$). The deeper reason for
the solvability may be rooted in the fact that
the symbolic dynamics of uncoupled Chebychev maps $T_N$ is very simple, consisting of
$N$ different symbols which all occur with the same probability and which are statisticaly
independent (similar to, but not exactly equal, to the observed symbolic dynamics of share price changes \cite{xu}). But while the symbols are statistically independent, the
actual iterates $x_n$ of the maps are not \cite{HC1,HC2}. Once again let us remark
that all Chebychev maps with $N \geq 2$ are conjugated to a Bernoulli shift and thus they are also
$\varphi$-mixing according to the definition given in section II.

Although the functional equations of the previous sections look quite complicated, it is remarkable that for certain mappings $T$ explicit solutions can be found.  These solutions
 are smooth and satisfy, in leading order, simple differential equations, describing
 the leading-order corrections to the Fokker-Flanck equation in the vicinity
 of the Gaussian limit case, which are due to deterministic chaotic effects as represented
 by the Perron-Frobenius equation.

 The Ulam map $T(x)=1-2x^{2}$ is the negative of the second-order Chebychev map.
Here we have $\langle x^{2} \rangle=1/2$ \cite{HC1}, \cite{HC2}, and the stationary solution of the Fokker-Planck equation $(41)$ is given by
\begin{equation}
	p_{0}(y,t)=\left(\frac{2}{\pi}\right)^{1/2}e^{-2y^{2}}
\end{equation}
This is to be expected, since the Ulam map satisfies a functional central limit theorem \cite{Tirn}.  To obtain the next-order correction term $\alpha(y,t)$, we first have to determine the inhomogeneous source term $(\partial/\partial y)\int dx xb(x,y,t)$ in Eq. $(42)$ by solving Eq. $(44)$.  Let us choose the following ansatz for the solutions of Eq. $(44)$:
\begin{equation}
	b(x,y,t)=h(x)\beta_{0}(y,t)+xh(x)\beta_{1}(y,t)
\end{equation}
Here $\beta_{0}$ and $\beta_{1}$ are appropriate functions independent of $x$.  Putting Eq. $(47)$ into Eq. $(44)$, we obtain on the left-hand side
\begin{equation}
\begin{aligned}
l ={} & h(x')\beta_{0}(y,t)+x'h(x')\beta_{1}(y,t) \\
	& = h(T(x))\left[\beta_{0}(y,t)+(1-2x^{2})\beta_{1}(y,t)\right]
\end{aligned}
\end{equation}
and on the right-hand side
\begin{equation}
\begin{aligned}
r ={} & \sum_{x\in T^{-1}(x')}\frac{1}{\vert T'(x)\vert}\biggl[h(x)\beta_{0}(y,t)+xh(x)\beta_{1}(y,t) \\
	& +h(x)\frac{1}{2}\left(x^{2}-\frac{1}{2}\right)\frac{\partial^{2}}{\partial y^{2}}p_{0}(y,t)\biggr] \\
	& =\sum_{x\in T^{-1}(x')}\frac{h(x)}{\vert T'(x)\vert}\left[\beta_{0}(y,t)-\frac{1}{4}(1-2x^{2})\frac{\partial^{2}}{\partial y^{2}}p_{0}(y,t)\right]
\end{aligned}
\end{equation}
From $l=r$ we get $\beta_{0}(y,t)$ arbitrary and
\begin{equation}
	\beta_{1}(y,t)=-\frac{1}{4}\frac{\partial^{2}}{\partial y^{2}}p_{0}(y,t)
\end{equation}
Thus, in the stationary case Eq. $(46)$ yields
\begin{equation}
	\beta_{1}(y,t)=\left(\frac{2}{\pi}\right)^{1/2}(1-4y^{2})e^{-2y^{2}}
\end{equation}

After a long calculation (details in Appendix A), one finds that all equations can be solved term by term and the final result for the stationary probability density of a stochastic differential equation driven by iterates of the Ulam map becomes
\begin{equation}
\begin{aligned}
	p(y) ={} & \left(\frac{2}{\pi}\right)^{1/2}\biggl[1+\tau^{1/2}\left(-\frac{8}{3}y^{3}+2y\right) \\
		& +\tau\left(\frac{32}{9}y^{6}-\frac{31}{3}y^{4}+\frac{15}{2}y^{2}-\frac{37}{48}\right)\biggr]e^{-2y^{2}}+O(\tau^{3/2})
\end{aligned}
\end{equation}
Our perturbative approach can (in principle) be extended to arbitrarily high orders in $\tau^{1/2}$.  It is interesting to notice that the non-Gaussian corrections of order $\tau^{k/2}, k=0,1,2$, multiplying the Gaussian function in Eq. (52) are of simple polynomial structure.  We may indeed use our approach to {\em define} an entire set of polynomials $\mathcal{P}_{k}$ by writing the $k$th-order correction term in the form $\tau^{k/2}(2/\pi)^{1/2}\mathcal{P}_{k}(y)$exp$(-2y^{2})$, similar as for wave functions in quantum mechanics.  The first three polynomials $\mathcal{P}_{0},\mathcal{P}_{1},\mathcal{P}_{2}$ are
\begin{equation}
	\mathcal{P}_{0}=1
\end{equation}
\begin{equation}
	\mathcal{P}_{1}=-\frac{8}{3}y^{3}+2y
\end{equation}
\begin{equation}
	\mathcal{P}_{2}=\frac{32}{9}y^{6}-\frac{31}{3}y^{4}+\frac{15}{2}y^{2}-\frac{37}{48}
\end{equation}
We would like to remark that the functional form given in Eq.~(52) has already been derived in \cite{HC2}; however, there it was obtained by a completely different graph-theoretic method.  The advantage of the general method presented here is that the Perron-Frobenius equation is applicable to an arbitrary dynamics driven by $T$, rather that just the special example $T(x)=1-2x^{2}$. Of course,
whether analytical solutions of the equations for other choices of $T$ and for more complicated
nonlinear relaxation dynamics can be found is another question; we found that the linear relaxation dynamics with
Chebychev maps is already complicated enough. Indeed we were able to find {\em analytic} solutions
for Chebychev maps only, but the equation can of course be solved numerically for other maps.

We once again notice \cite{Beck3,Beck4,Beck5} that the Ulam map is distinguished in comparison to other mappings.  For this map the first-order correction term $\alpha(y,t)$ obeys a relatively simple differential equation, namely the inhomogeneous Fokker-Planck equation
\begin{equation}
	\frac{\partial}{\partial y}(y\alpha(y,t))+\frac{1}{4}\frac{\partial^{2}}{\partial y^{2}}\alpha(y,t)-\frac{\partial}{\partial t}\alpha(y,t)=-\frac{1}{8}\frac{\partial^{3}}{\partial y^{3}}p_{0}(y,t)
\end{equation}
Moreover, the second-order correction term $\beta(y,t)$ obeys the relatively simple equation $(83)$ presented in Appendix A.  In general, for some arbitrary mapping $T$ the correction terms obey the much more complicated equations $(27)-(29)$.

There is a minor sign difference between the stationary probability density for the Ulam map and for the second-order Chebyshev polynomial.  Due to the Ulam map being the negative of the second-order Chebyshev polynomial the difference will show in the perturbative, first order correction term $\alpha(y)$ only - it
has opposite sign for both maps, i.e., the result for the stationary probability density of $T_{2}$ is
\begin{equation}
\begin{aligned}
	p(y) ={} & \left(\frac{2}{\pi}\right)^{1/2}\biggl[1+\tau^{1/2}\left(\frac{8}{3}y^{3}-2y\right) \\
		& +\tau\left(\frac{32}{9}y^{6}-\frac{31}{3}y^{4}+\frac{15}{2}y^{2}-\frac{37}{48}\right)\biggr]e^{-2y^{2}}+O(\tau^{3/2})
\end{aligned}
\end{equation}
\section{Explicit Solution for $N=3$, $T_{3}(x)=4x^{3}-3x$}
We will now look at the explicit solution for the $3^{rd}$ order Chebyshev polynomial $T(x)=4x^{3}-3x$.  As $h(x)$ is the same for all Chebyshev polymials \cite{Beck6}, \cite{HC1}-\cite{Tirn}, we have $\langle x^{2}\rangle=1/2$ and $\langle x^{4}\rangle=3/8$, and the stationary solution of the Fokker-Planck equation $(41)$ is a Gaussian with the same variance for all Chebyshev maps.  Chebyshev maps $T_N$ with $N$ odd
are odd functions satisfying $T_N(-x)=-T_N(x)$, and as a consequence of this symmetry
all velocity distributions $p(y)$ are symmetric in $y$.
Again all functional equations of our perturbative approach can be solved, see Appendix B for details.

As a consequence of Lemmas 1-3 (see Appendix C), any Chebyshev map with $N\geq 2$ satisfies
\begin{equation} \nonumber
        \varphi(x',y,t)=\sum_{x\in T^{-1}(x')} \frac{1}{|T'(x)|}\varphi(x,y,t)
\end{equation}
which is solved by $\varphi(x,y,t)=h(x)p_{0}(y,t)$. The equation
\begin{equation} \nonumber
        a(x',y,t)=\sum_{x\in T^{-1}(x')} \frac{1}{|T'(x)|}a(x,y,t)
\end{equation}
is solved by $a(x,y,t)=h(x)\alpha(y,t)$; and the equation
\begin{equation}
        b(x',y,t)=\sum_{x\in T^{-1}(x')} \frac{1}{|T'(x)|}\left[b(x,y,t)+\frac{1}{4}h(x)\left(2x^{2}-1\right)\frac{\partial^{2}}{\partial y^{2}}p_{0}\left(y,t\right)\right]
\end{equation}
which is derived from $(17)$ and $(32)$ will be discussed in the appendix (Lemma 4).

The final result of a long calculation (details in Appendix B), is that in leading order the stationary probability density of a particle kicked by $T_3$ is given by
\begin{equation}
	p(y)=\left(\frac{2}{\pi}\right)^{1/2}\left(1+\tau\left(\frac{1}{3}y^{4}+\frac{3}{2}y^{2}-\frac{7}{16}\right)\right)e^{-2y^{2}} +O(\tau^{2}).
\end{equation}
\section{Explicit solution for $N=4$, $T_{4}(x)=8x^{4}-8x^{2}+1$}
      Let us now determine the explicit solution for the $4^{th}$ order Chebyshev polynomial $T(x)=8x^{4}-8x^{2}+1$.  Again, as in the last section, we have the same values for $h(x)$,$\langle x^{2}\rangle$ and $\langle x^{4}\rangle$ \cite{Beck6}, \cite{HC1}-\cite{Tirn}; and the same stationary solution of the Fokker-Planck equation $(41)$.

For this $T$, our first-order correction term satisfies $\alpha(y,t)=0$, so we will focus on finding $\beta(y,t)$.  We have
\begin{equation} \nonumber
	\frac{\partial}{\partial y} \int dx \ xc(x,y,t)=0
\end{equation}
\begin{equation} \nonumber
	\frac{1}{2}\frac{\partial^{2}}{\partial y^{2}} \int dx \ x^{2}b(x,y,t)=\frac{1}{2}\langle x^{2}\rangle \frac{\partial^{2}}{\partial y^{2}}\beta(y,t)=\frac{1}{4}\frac{\partial^{2}}{\partial y^{2}}\beta(y,t)
\end{equation}
\begin{equation} \nonumber
	\frac{\partial}{\partial t}\beta(y,t)=0
\end{equation}
\begin{equation} \nonumber
	\frac{\partial}{\partial t}p_{0}(y,t)-\frac{1}{2}\frac{\partial^{2}}{\partial t^{2}}p_{0}(y,t)=0
\end{equation}
Therefore, $(43)$ becomes
\begin{equation}
	0=\frac{\partial}{\partial y}(y\beta(y,t))+\frac{1}{4}\frac{\partial^{2}}{\partial y^{2}}\beta(y,t)+\biggl[\frac{1}{2}y\frac{\partial}{\partial y}+\frac{1}{2}(y^{2}+1)\frac{\partial^{2}}{\partial y^{2}}+\frac{1}{4} y\frac{\partial^{3}}{\partial y^{3}}+\frac{1}{64}\frac{\partial^{4}}{\partial y^{4}}-\frac{1}{2}\biggr]p_{0}(y,t)
\end{equation}
which, in the stationary case becomes
\begin{equation}
	0=\frac{\partial}{\partial y}(y\beta(y))+\frac{1}{4}\frac{\partial^{2}}{\partial y^{2}}\beta(y)+\left(-4y^{4}+10y^{2}-\frac{7}{4}\right)\left(\frac{2}{\pi}\right)^{1/2}e^{-2y^{2}}.
\end{equation}
One can easily check that this differential equation is solved by
\begin{equation}
	\beta(y)=\left(\frac{2}{\pi}\right)^{1/2}\left(-y^{4}+\frac{7y^{2}}{2}+C\right)e^{-2y^{2}}
\end{equation}
where $C$ is a constant. One must have $\int dy \ \beta(y)=0$ and from this one arrives at

\begin{equation}
	\left(-\frac{3}{16}+\frac{7}{2} \cdot \frac{1}{4} +C\right) = 0 \Rightarrow C=-\frac{11}{16}
\end{equation}
The final result for the stationary probability density is thus
\begin{equation}
	p(y)=\left(\frac{2}{\pi}\right)^{1/2}\left(1+\tau\left(-y^{4}+\frac{7y^{2}}{2}-\frac{11}{16}\right)\right)e^{-2y^{2}}+O(\tau^{3/2})
\end{equation}

The same result applies to any higher Chebychev map with $N>4$,
as the only differences occuring in this case are within the higher-order terms represented by $O(\tau^{3/2})$ or higher.
\section{Numerical experiments}
We will now compare our analytically derived results with numerically generated histograms, obtained by iterating the map (1). For all simulations considered in this paper $x_{0}$ is chosen randomly with
uniform distribution such that $x_{0}\in(-1,1)$, and $y_{0}=0$.  All of our iterations were programmed in Matlab, with all histograms built from the $y$-values in Eq.(1) using $10^{3}$ bins and normalized to 1.
		
		First we look at uncoupled Chebyshev maps of order $N=2,3,4$ acting as a deterministic chaotic driving force.  In our numerical experiment the map (1) is iterated $3\times10^{7}$ times.
          In Figs 1a-1f, we see that our analytic results, Eqs.$(57)$, $(59)$, and $(64)$, tightly fit the numerical histograms.  We can see  that the smaller $\tau$, the better the coincidence. Figures 1a and 1b show numerical histograms $\rho (y)$ (in the physical
          picture, these are particle velocity distributions) for driving forces generated by $T_{2}$, Figures 1c and 1d show the same for $T_{3}$, and Figures 1e and 1f for $T_{4}$.  Note, any disparity between the histograms and the analytic prediction comes from the neglection of the higher order terms of $O(\tau^{3/2})$ in the fits.

\newpage

\vspace{-6cm}

\begin{figure}[h]
        \subfloat[]{\includegraphics[width=0.4075\textwidth]{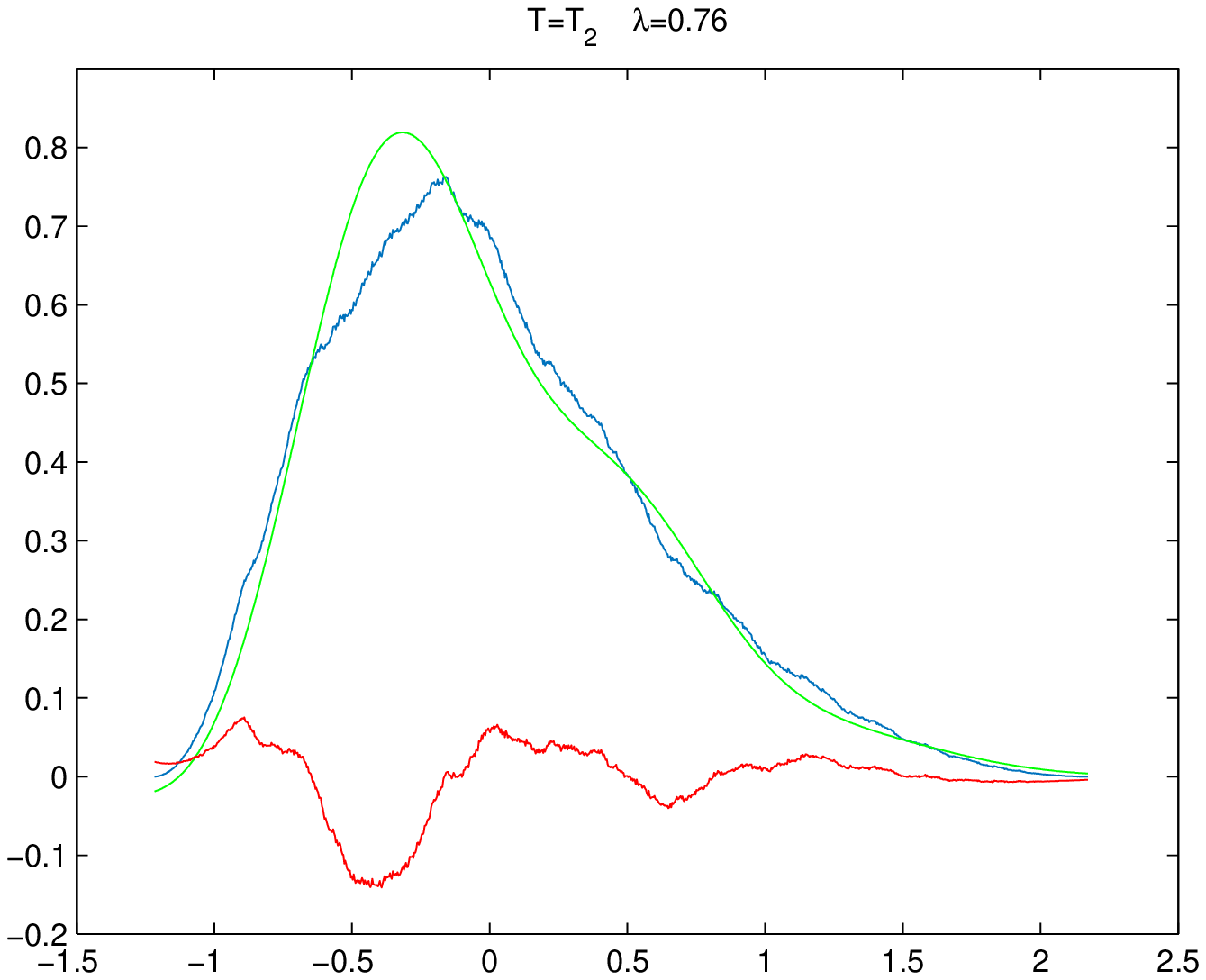}}
        \subfloat[]{\includegraphics[width=0.4075\textwidth]{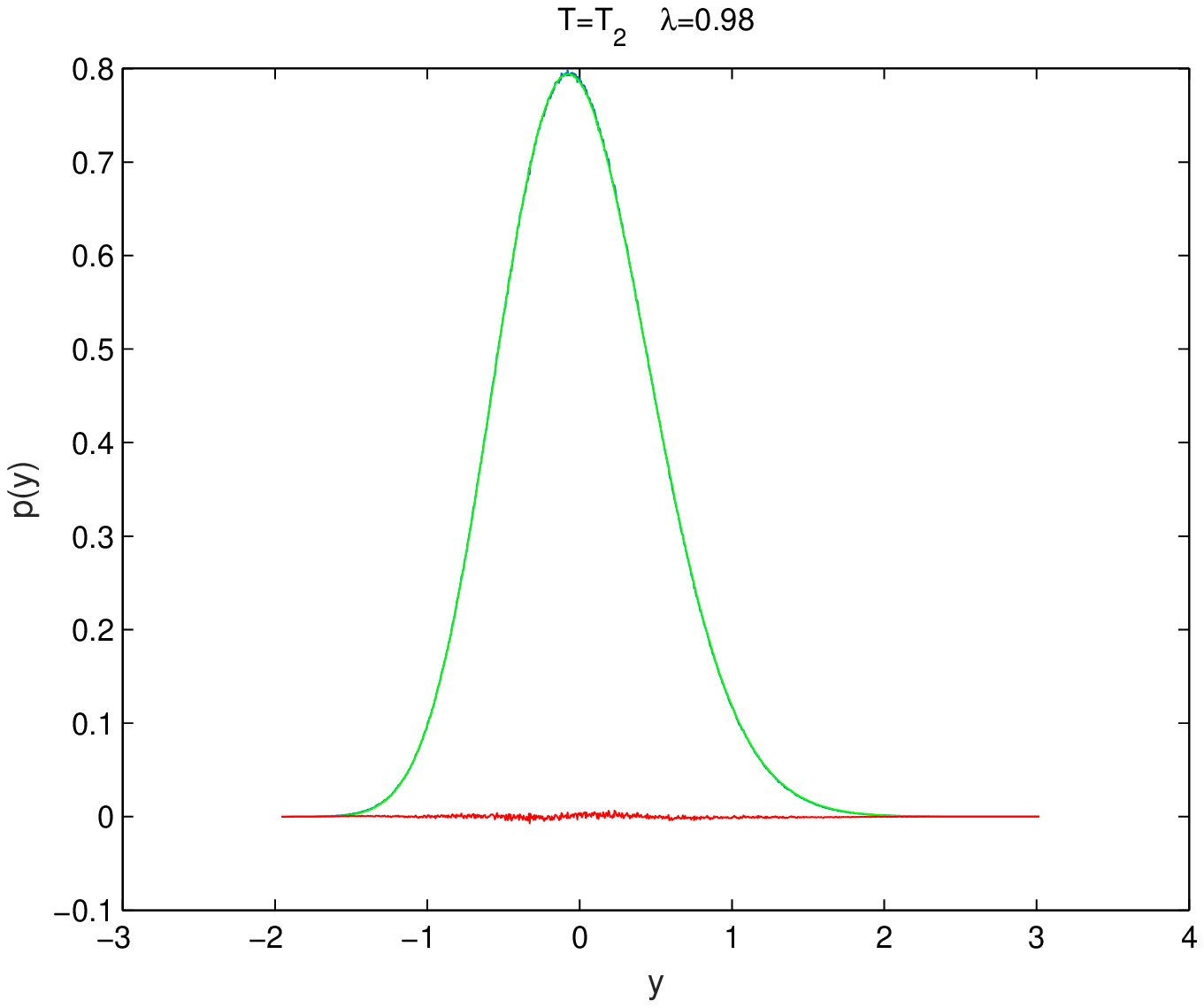}}\\
        \subfloat[]{\includegraphics[width=0.4075\textwidth]{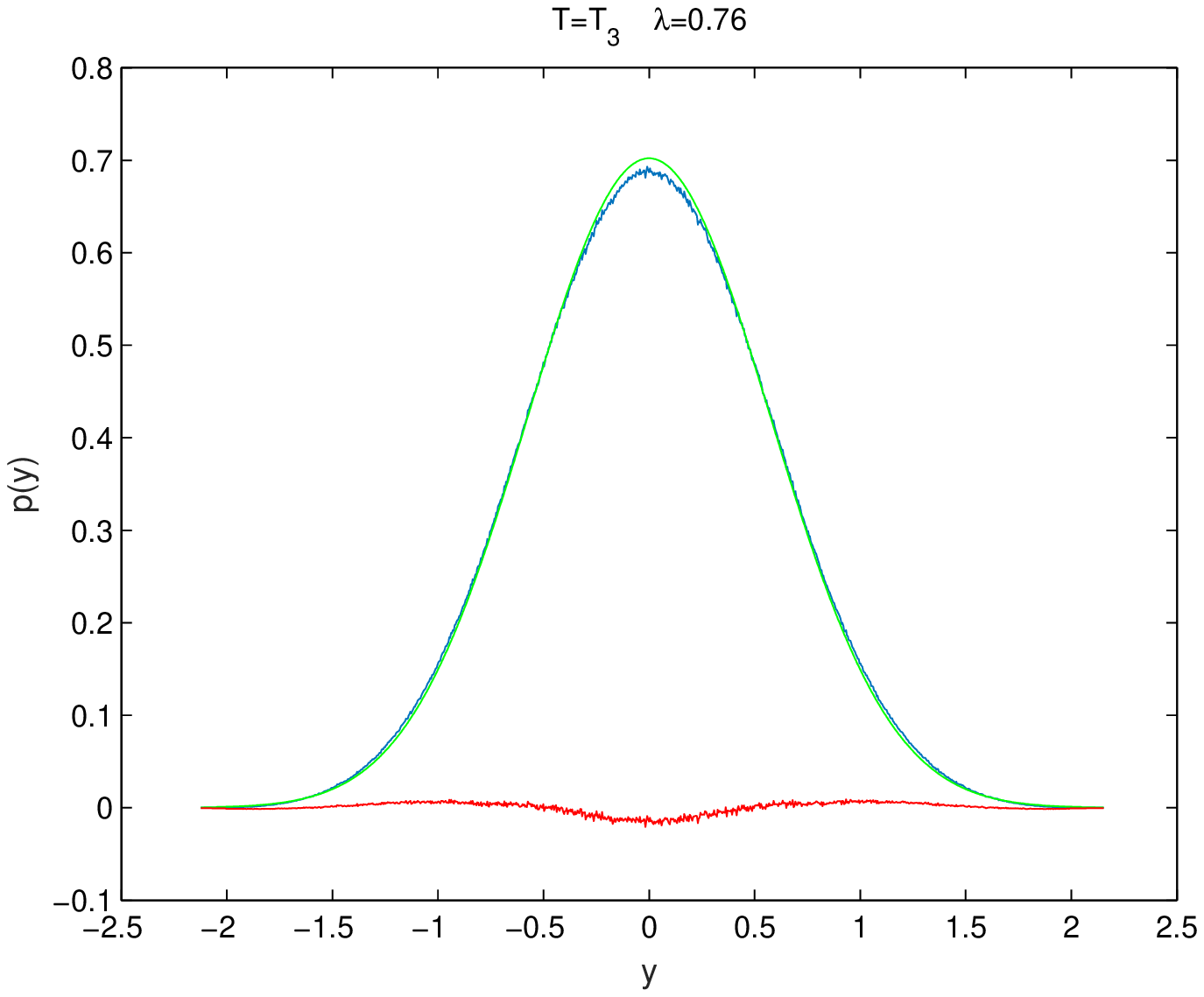}}
        \subfloat[]{\includegraphics[width=0.4075\textwidth]{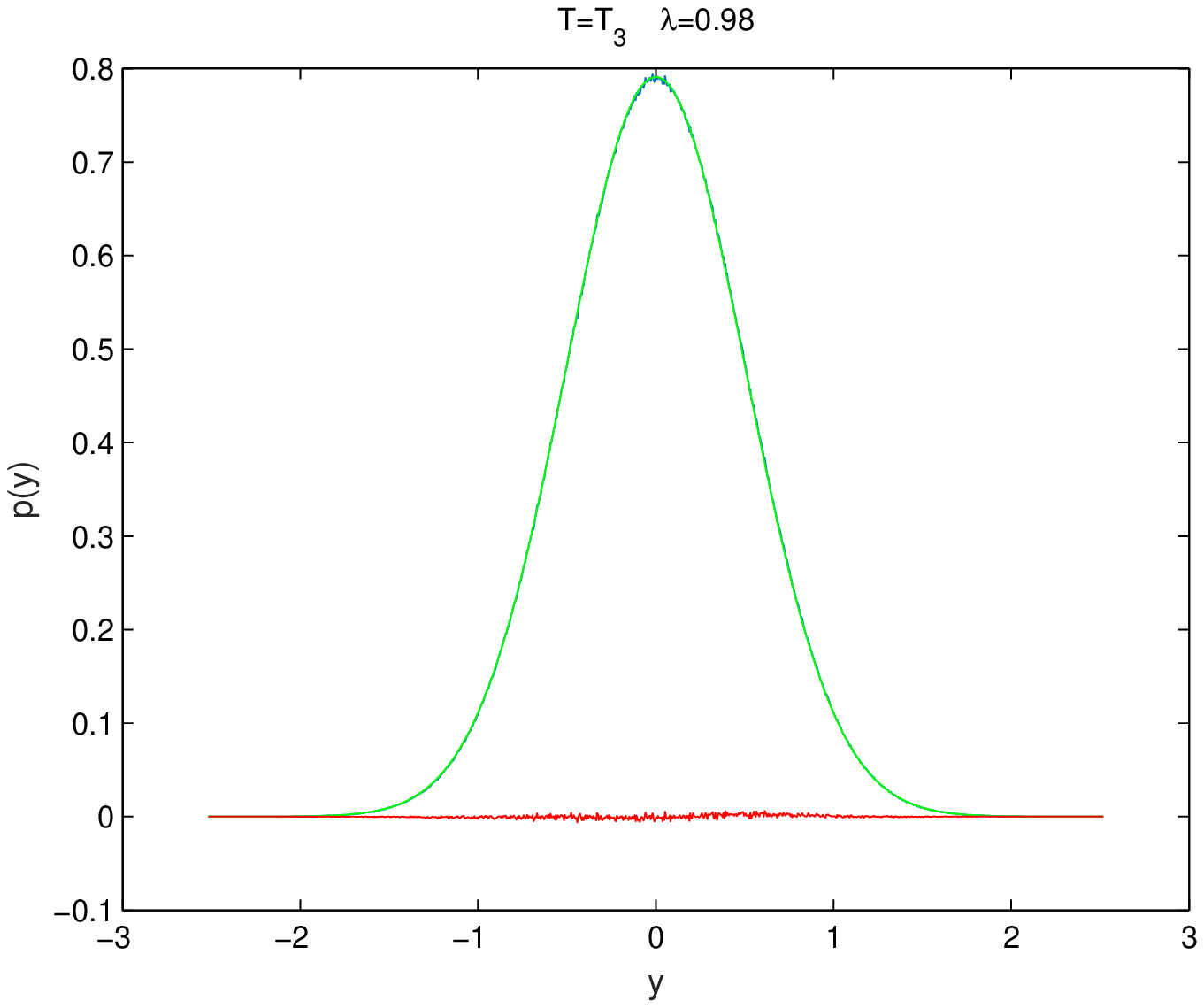}}\\
        \subfloat[]{\includegraphics[width=0.4075\textwidth]{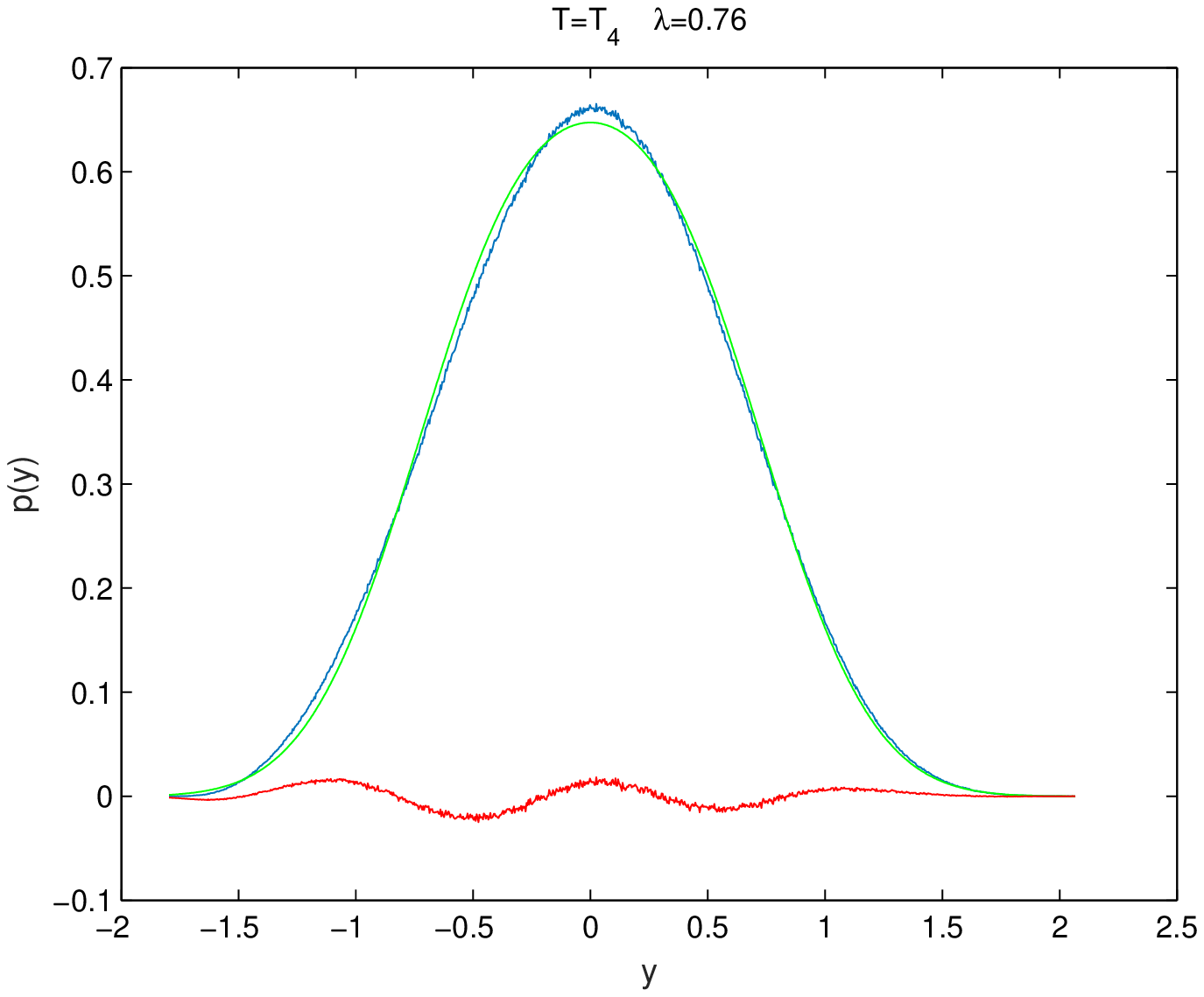}}
        \subfloat[]{\includegraphics[width=0.4075\textwidth]{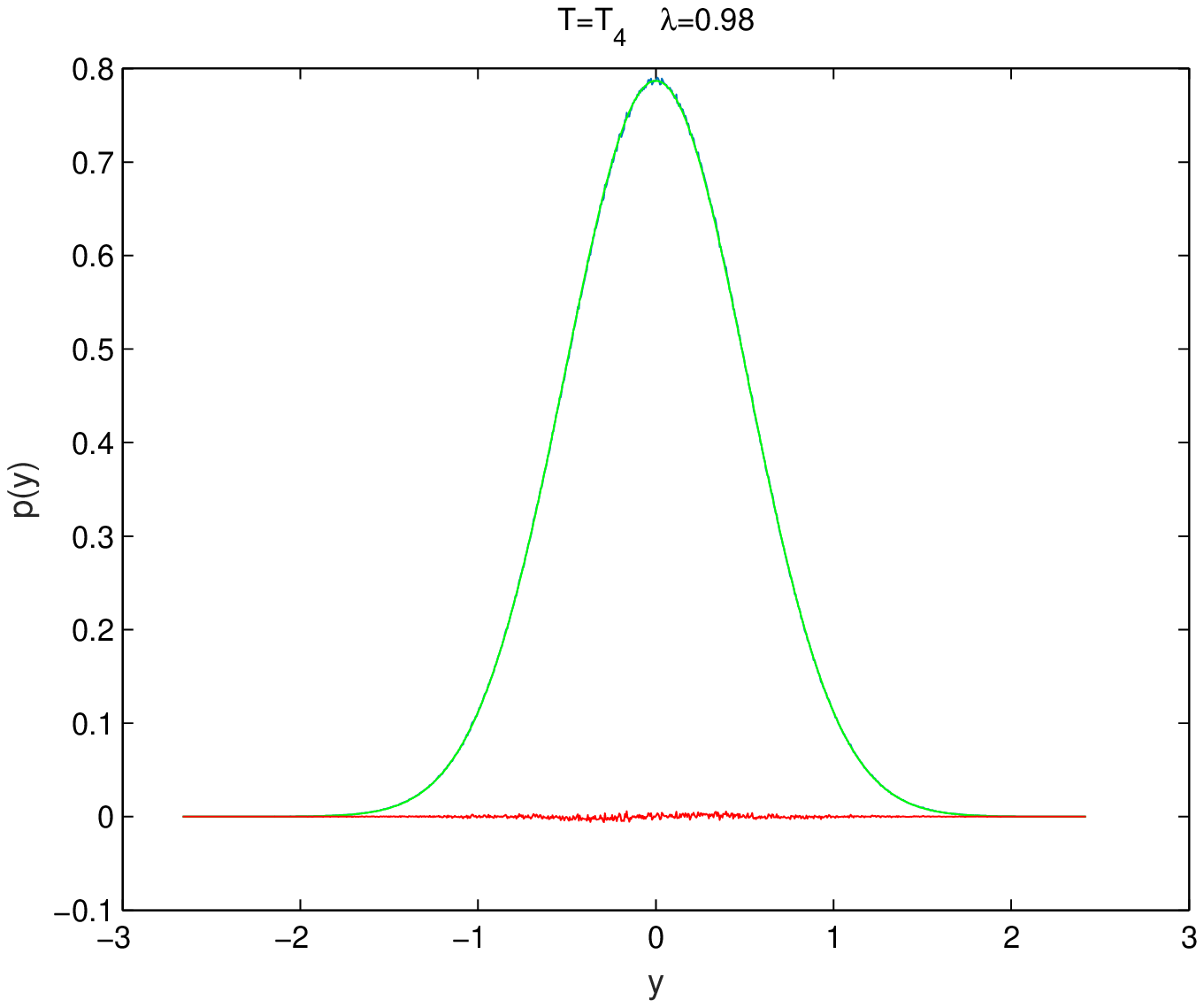}}

\vspace{4cm}

        \caption{Comparison of the numerical results (blue) with the analytic results (green) for $T_{2}$ (a),(b); $T_{3}$ (c),(d); $T_{4}$ (e),(f). The subtractions between these results, denoting the tightness of fit, are shown in red.  The results shown are for $\lambda=e^{-\tau}=0.76$ (a),(c),(e); and $\lambda=e^{-\tau}=0.98$ (b),(d),(f).}

\end{figure}

\newpage
         The figure shows the histograms and the theoretical prediction for two values $\lambda=0.76, 0.98$ so as to illustrate the behavior of the system when $\lambda$ is both further away and closer to 1 ($\tau$ is further away and closer to 0).  Clearly, for $\lambda\rightarrow 1$ there is convergence to a Gaussian, due to the $\varphi$-mixing property
         mentioned earlier.
    We also show the differences between the histograms and our analytic curves, shown in the figure in red.  It can be seen that as $\lambda\rightarrow 1$, for all three Chebyshev polynomials the difference becomes smaller and smaller,
consistent with convergence to zero.

Another interesting aspect of our analytic work in the previous sections was the extraction of the different polynomials that
describe the leading order corrections when approaching the Gaussian limit case. These polynomials
can be regarded as eigenfunctions of a rescaling operator $\tau \to \frac{1}{2} \tau$ that describes the {\em approach} to the Gaussian limit case
under rescaling of the parameter $\tau$, in the vicinity of the Gaussian fixed point.
For $T_2$ we obtain after re-arranging
               \begin{equation}
	\left(\frac{\sqrt{\frac{\pi}{2}}\rho(y)}{e^{-2y^{2}}}-1\right)\frac{1}{\sqrt{\tau}}=\frac{8}{3}y^{3}-2y+\sqrt{\tau}\left(\frac{32}{9}y^{6}-\frac{31}{3}y^{4}+\frac{15}{2}y^{2}-\frac{37}{48}\right)
+O(\tau)\end{equation}
so if we choose for $\rho (y)$ a numerically generated histogram one can extract
the polynomials by the above re-arrangement. This is shown in Figs 2a and 2b, taking for $\rho(y)$ the numerically obtained histogram and two different values of $\tau$. Excellent agreement between numerical results and our
analytic predictions is obtained for $\tau$ small.
Similarly, we obtain for $T_3$
\begin{equation}
	\left(\frac{\sqrt{\frac{\pi}{2}}\rho(y)}{e^{-2y^{2}}}-1\right)\frac{1}{\tau}=\frac{1}{3}y^{4}+\frac{3}{2}y^{2}-\frac{7}{16}
+O(\tau)
\end{equation}
and for $T_4$
\begin{equation}
	\left(\frac{\sqrt{\frac{\pi}{2}}\rho(y)}{e^{-2y^{2}}}-1\right)\frac{1}{\tau}=-y^{4}+\frac{7y^{2}}{2}-\frac{11}{16}
+O(\tau^{1/2})\end{equation}

Figs 2a, 2b confirm the analytically calculated polynomial form of the leading correction terms
for $T_2,T_3,T_4$, and as expected the deviations
from the calculated polynomials (the higher order terms) become very small for small $\tau$.  For these simulations, we iterated the map (1) $1.75\times10^{7}$ times.

\begin{figure}
    \centering
    \subfloat[]{\includegraphics[width=.5\linewidth]{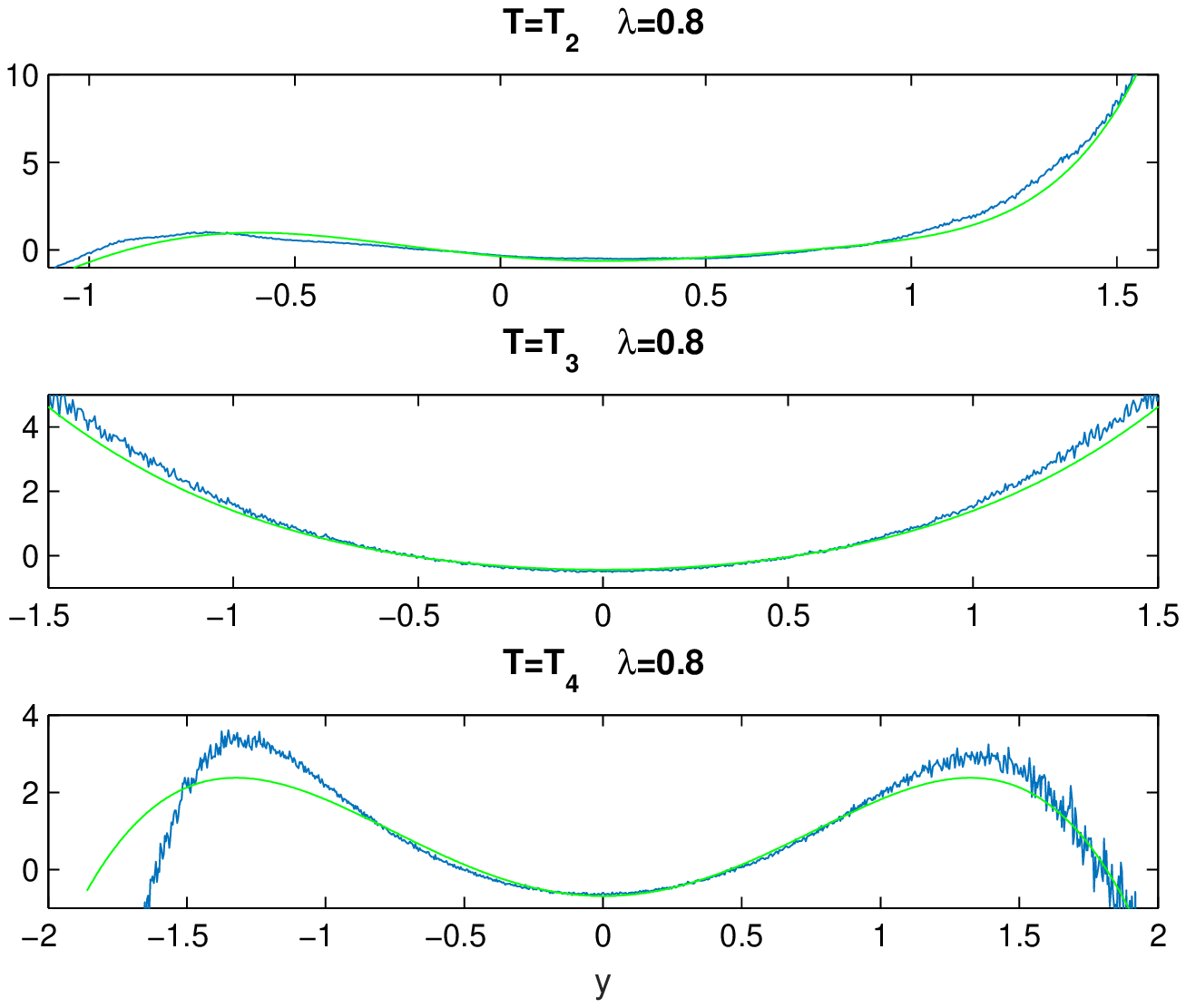}\label{fig:loop}} \hfill
    \subfloat[]{\includegraphics[width=.5\linewidth]{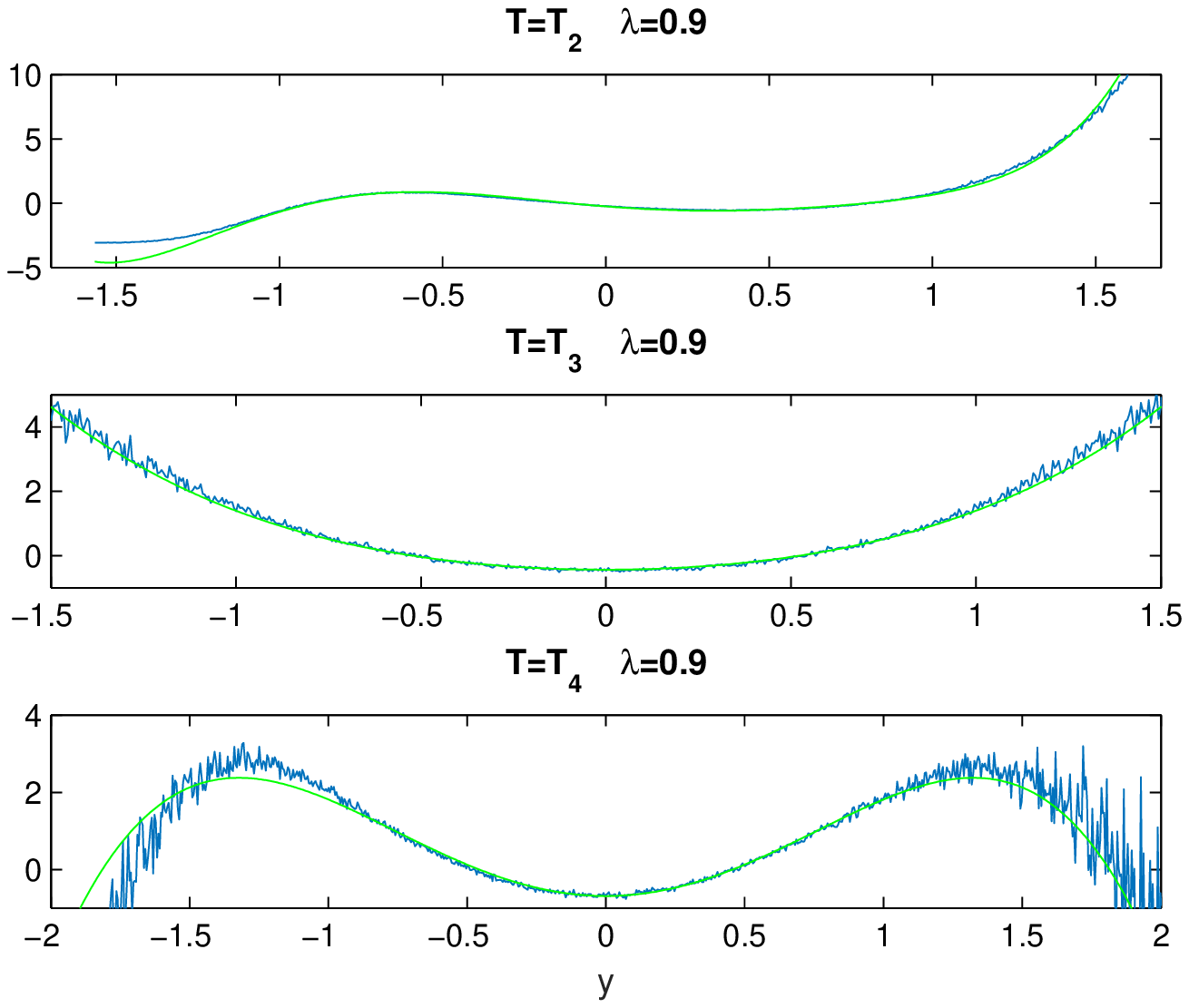}\label{fig:twopath}}

    \vspace{7cm}

    \caption{Leading-order corrections to the Gaussian limit case for $\lambda=e^{-\tau}=0.8$ (a), and $\lambda=e^{-\tau}=0.9$ (b).  The numerical results are shown in blue and the calculated polynomials are shown in green.  The Chebyshev map that each plot relates to is demarcated in the plot.}
    \label{fig:routing}
\end{figure}

Finally, we may ask what happens if the Chebyshev maps that drive the system are coupled. Here we use
diffusively coupled map lattices of the form
\begin{equation}
x_{n+1}^{(i)}=(1-\alpha) T_N(x_n^{(i)}) + \frac{1}{2} \alpha (T_N (x_n^{(i-1)})+T_N(x_n^{(i+1)})),
\end{equation}
with periodic boundary conditions and lattices of size 1000.
For large couplings $\alpha$ these systems may lose the $\varphi$-mixing property, hence convergence to
a Gaussian process is not guaranteed anymore if components of this coupled map lattice drive the
linear relaxation dynamics. In fact for finite $\tau$ many complicated shapes of densities
can be created. To quantitatively measure the deviation from a Gaussian, Figs 3a-3c show the kurtosis $\kappa$
of the numerically obtained histograms $\rho(y)$ as a function of the coupling $\alpha$, with $\alpha\in[0,1]$, and $\Delta \alpha=.005$, for various values of $\tau$. All lattice sites were used as driving forces, iterated $10^{5}$ times giving $10^{8}$ data points, and all the driven velocities were put into the
same histogram.  This was done for $\lambda=0.6,0.7,0.8,0.9$.
Some example densities are shown in Figs 4a-4f.  Further, the log-plots of these densities are also shown in Figs 5a-5f.  These were simulated with $2\times10^{5}$ iterations giving $2\times10^{8}$ data points for $\alpha=0.0,0.5,1.0$ and $\lambda=0.6,0.9$.  These plots nicely illustrate  the deviation from a Gaussian as the semi-logarithmic plot of a Gaussian is a parabola.  Convergence to the Gaussian limit case (with has $\kappa=3$)
is mathematically rigorously proven for $\alpha=0$ only. The other cases of $\alpha$ can
  lead to more complex non-Gaussian behavior due to strong temporal correlations of the iterates,
  where the ordinary central limit behavior is not guaranteed anymore \cite{Tirn}.
  An interesting example is the case of $T_3$ with $\alpha =0.5$ where our numerics
  indicates that complex non-Gaussian behavior seems to persist even for $\tau \to 0$.
  In the quantum field theoretical applications
described in \cite{Beck7}-\cite{SM} this corresponds to a strongly interacting field theory rather than a free field.

\begin{figure}[h]
        \subfloat[]{\includegraphics[width=0.475\textwidth]{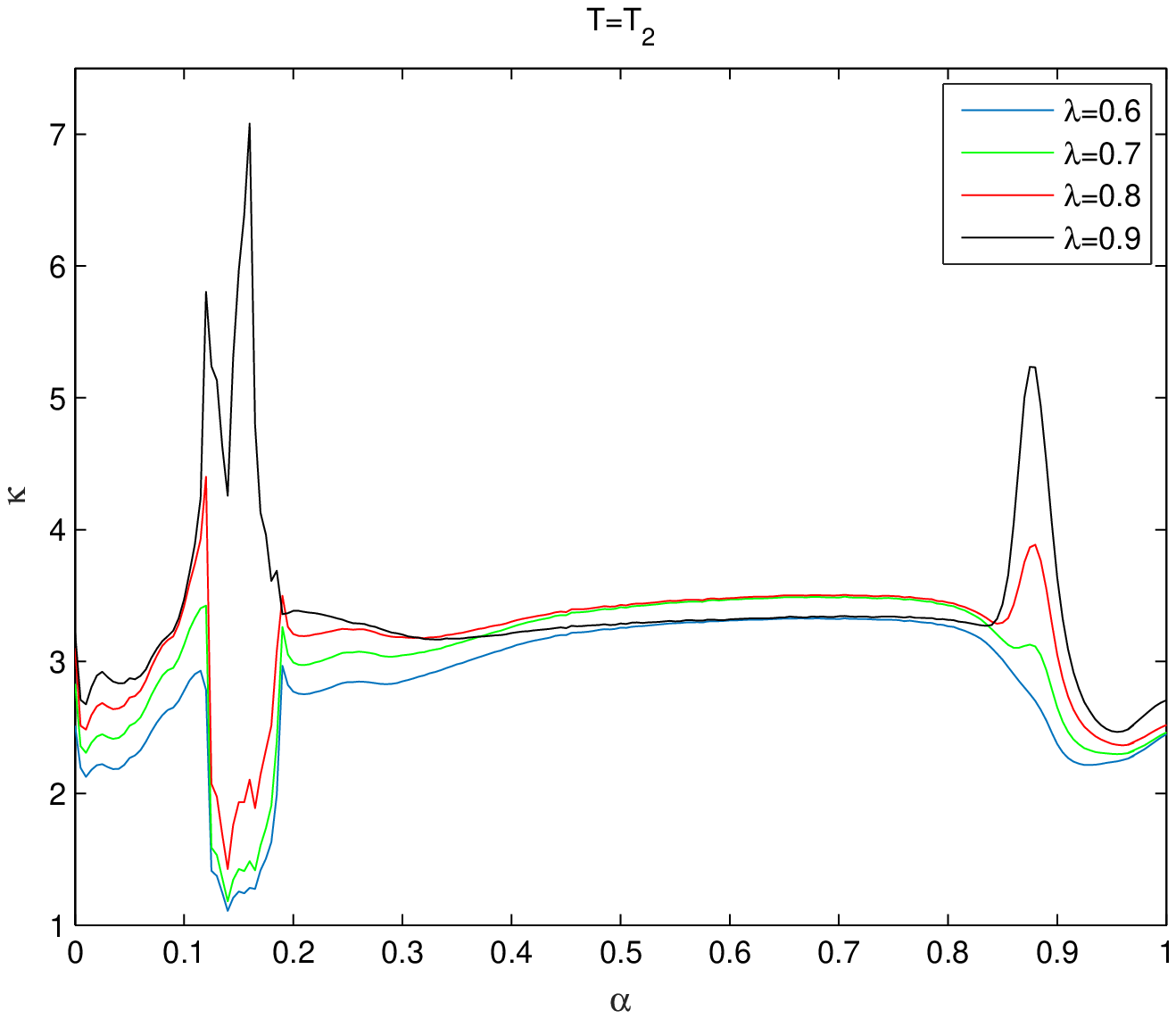}}
        \subfloat[]{\includegraphics[width=0.475\textwidth]{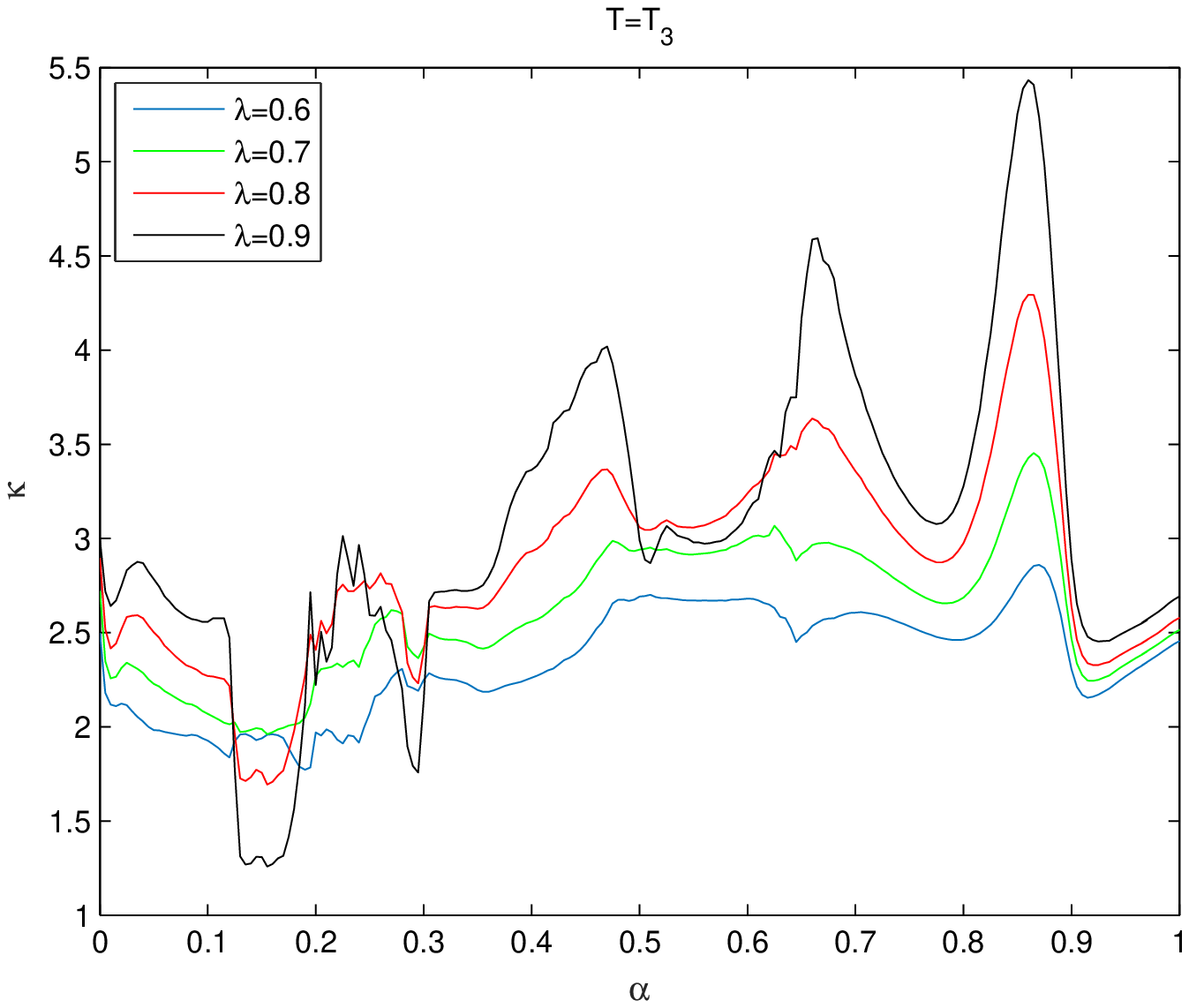}}\\
        \subfloat[]{\includegraphics[width=0.475\textwidth]{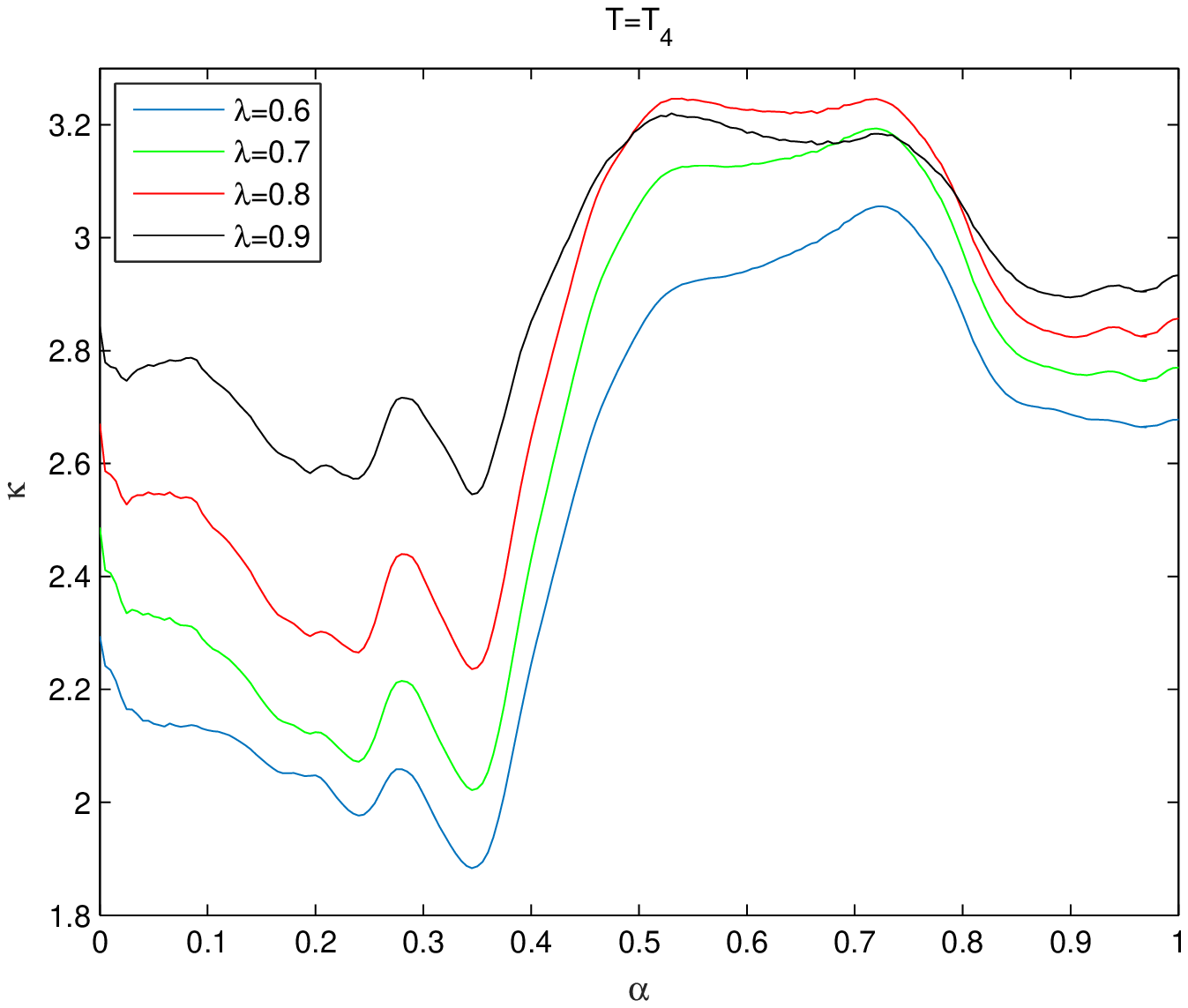}}

        \vspace{5cm}

        \caption{Kurtosis $\kappa$ of the velocity distributons as a function of the couplings $\alpha$. Each plot corresponds to a different Chebyshev polynomial being coupled: $T_{2}$ (a), $T_{3}$ (b), $T_{4}$ (c).  The blue lines describe $\lambda=0.6$, green describes $\lambda=0.7$, red describes $\lambda=0.8$, and black describes $\lambda=0.9$.  $\Delta \alpha=.005$. For some parameter values $\alpha$, there are particularly strong deviations of the kurtosis
        from the Gaussian value $\kappa =3$.}
\end{figure}
\begin{figure}[h]
        \subfloat[]{\includegraphics[width=0.4075\textwidth]{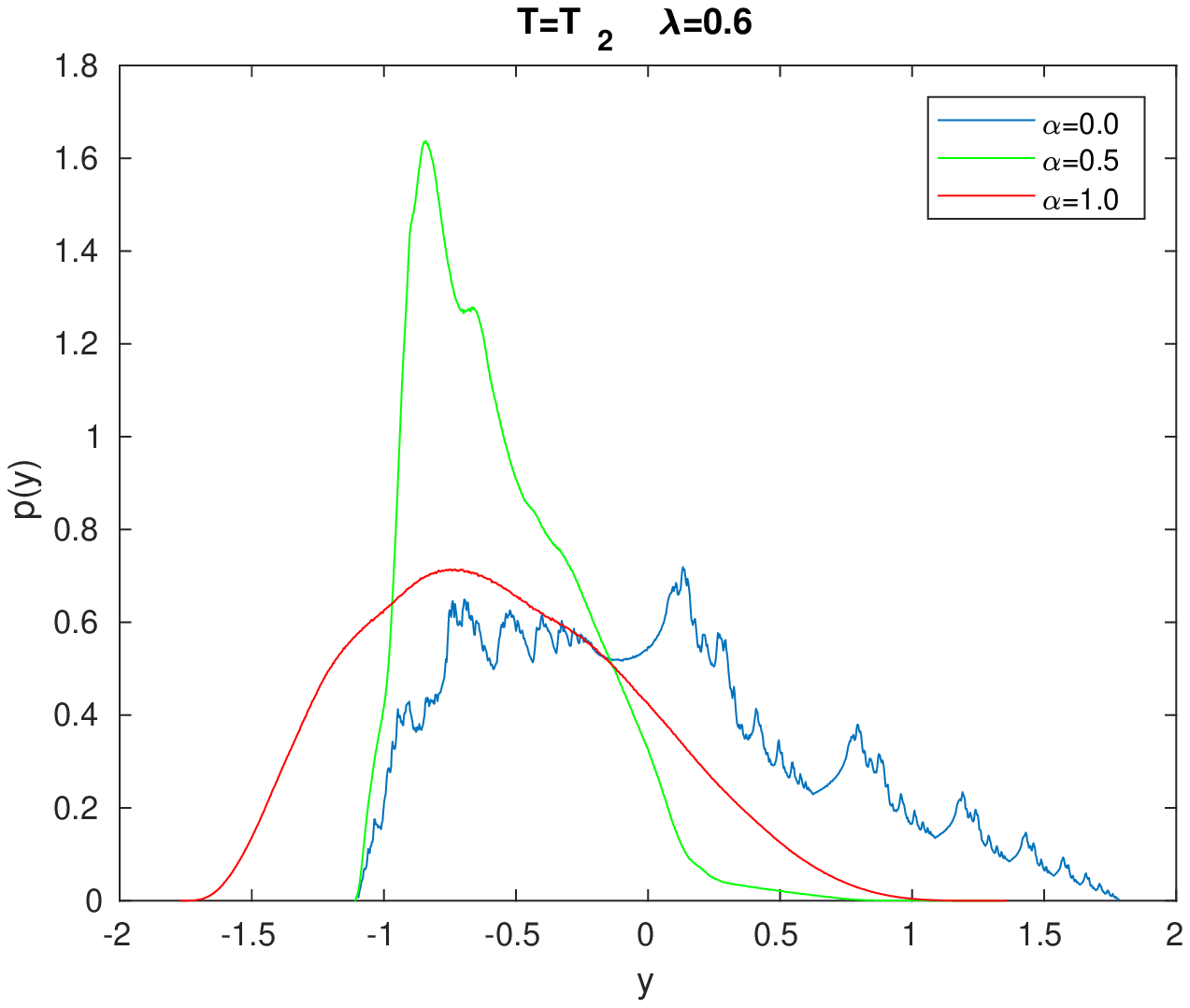}}
        \subfloat[]{\includegraphics[width=0.4075\textwidth]{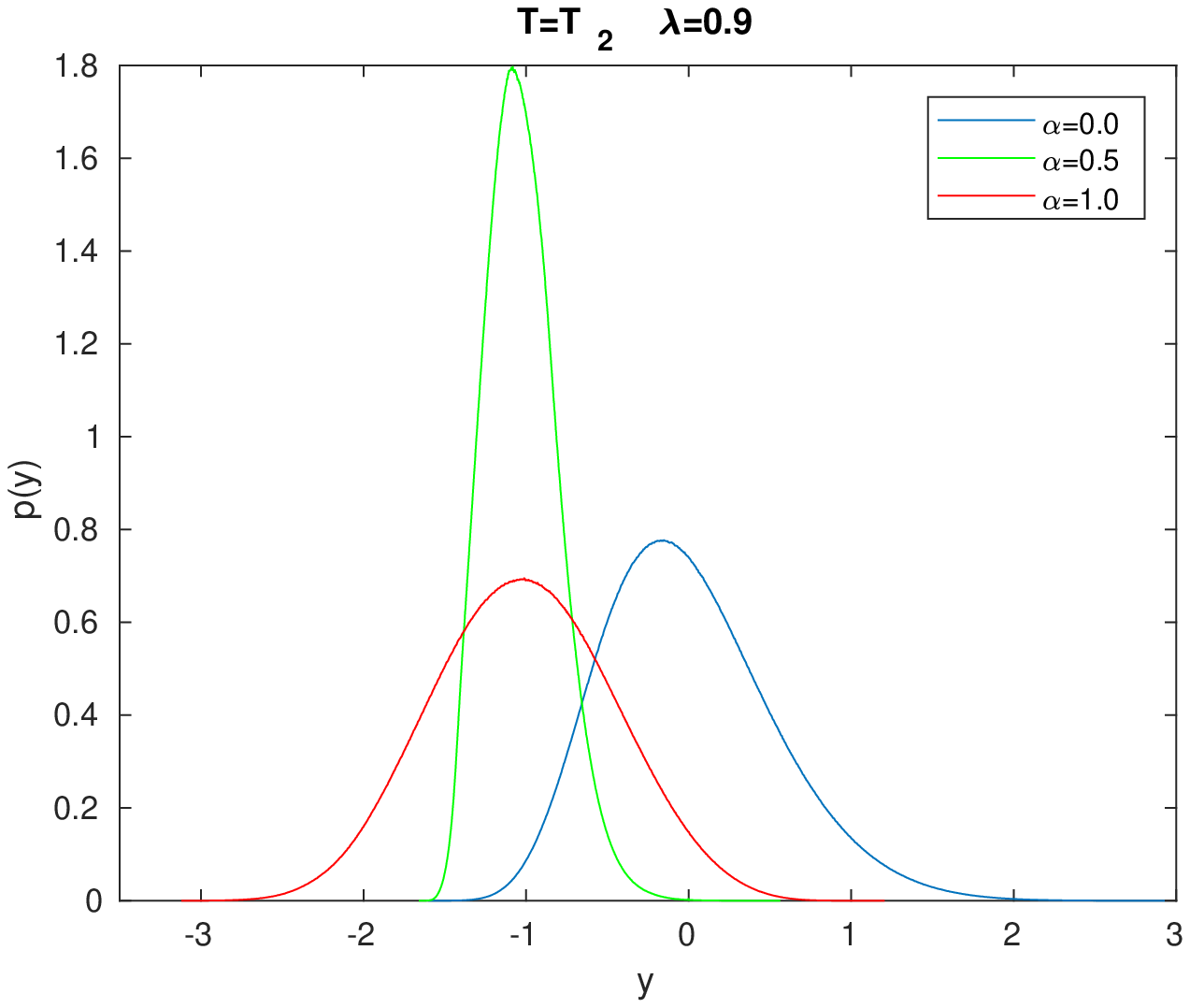}}\\
        \subfloat[]{\includegraphics[width=0.4075\textwidth]{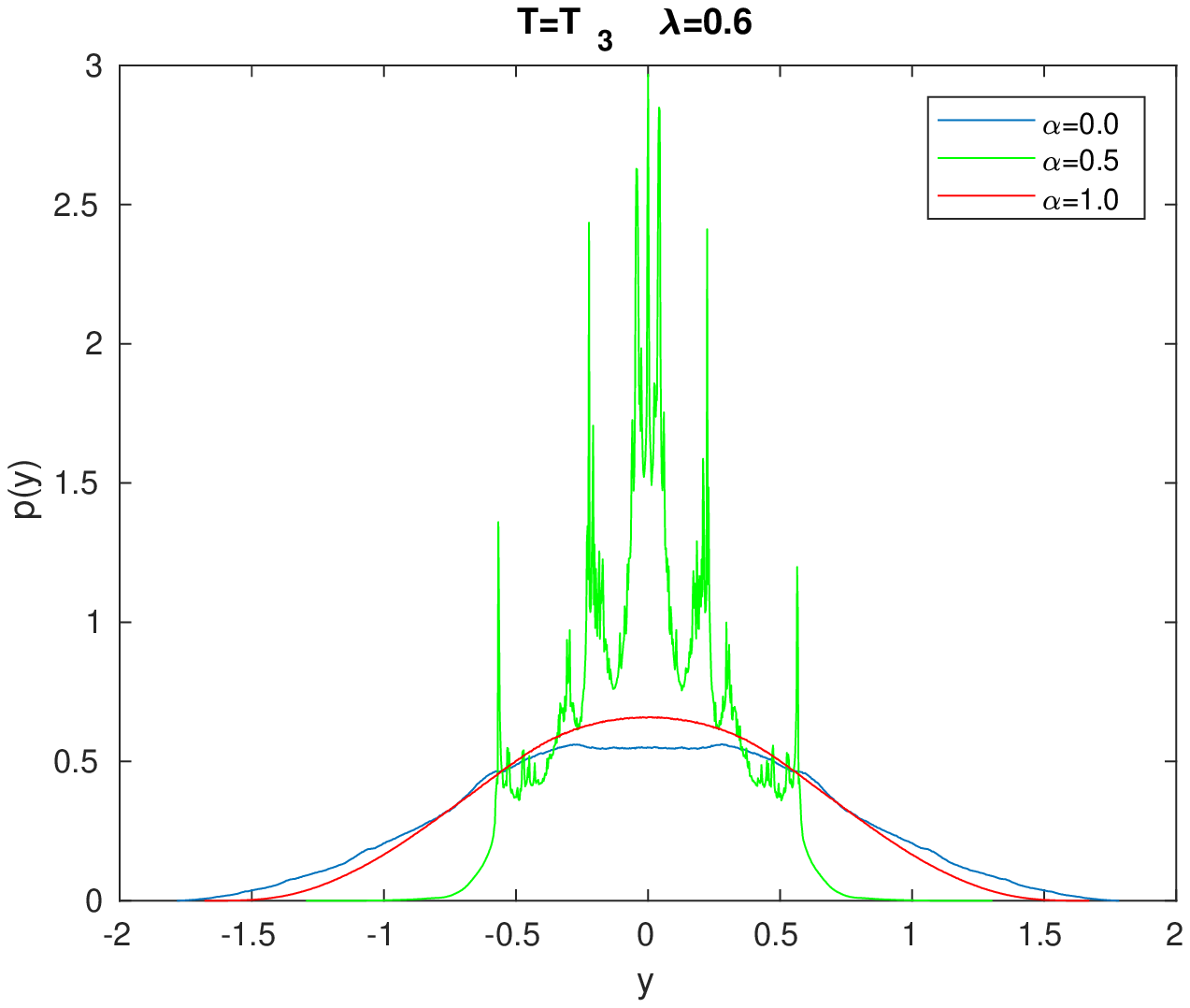}}
        \subfloat[]{\includegraphics[width=0.4075\textwidth]{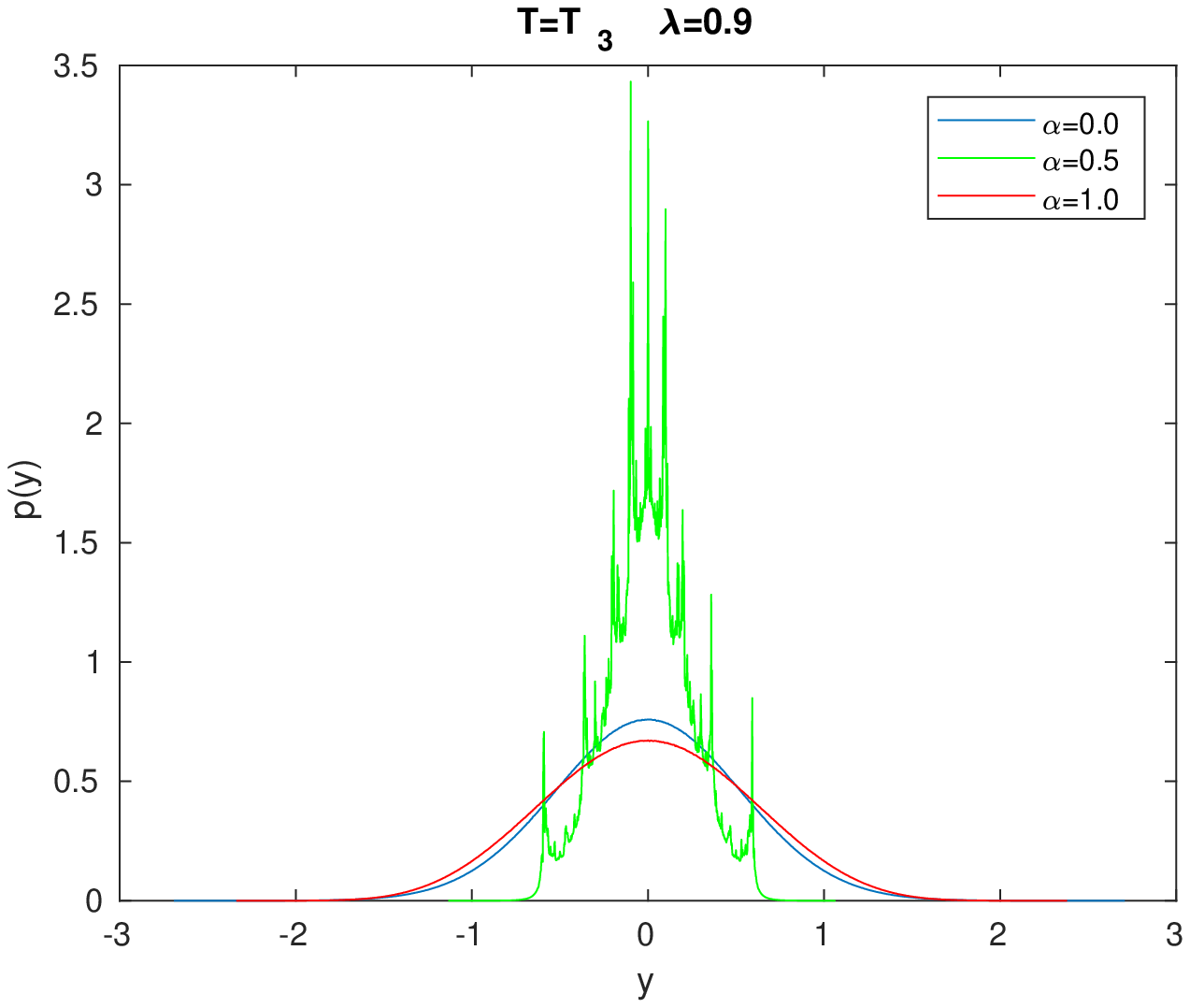}}\\
        \subfloat[]{\includegraphics[width=0.4075\textwidth]{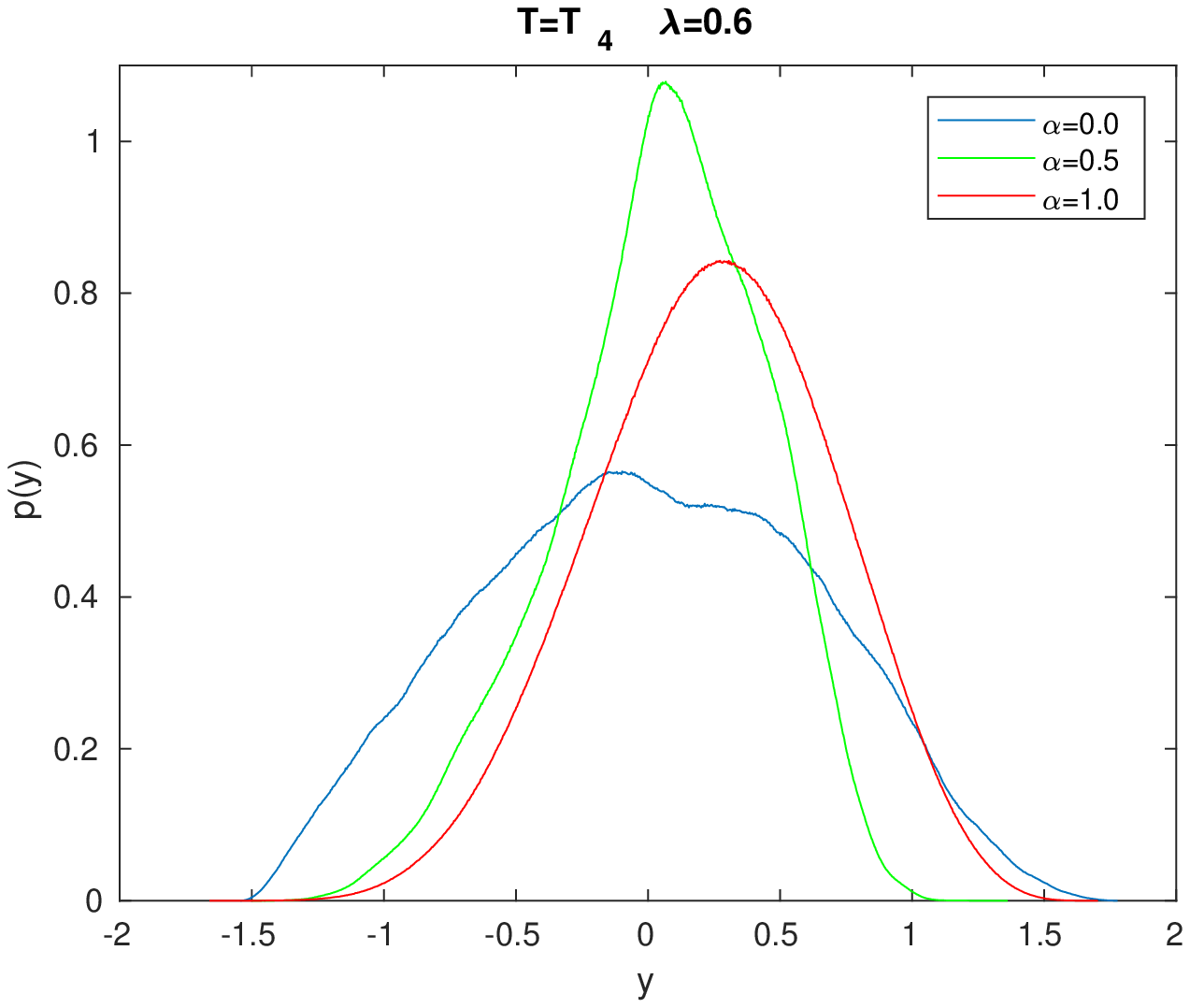}}
        \subfloat[]{\includegraphics[width=0.4075\textwidth]{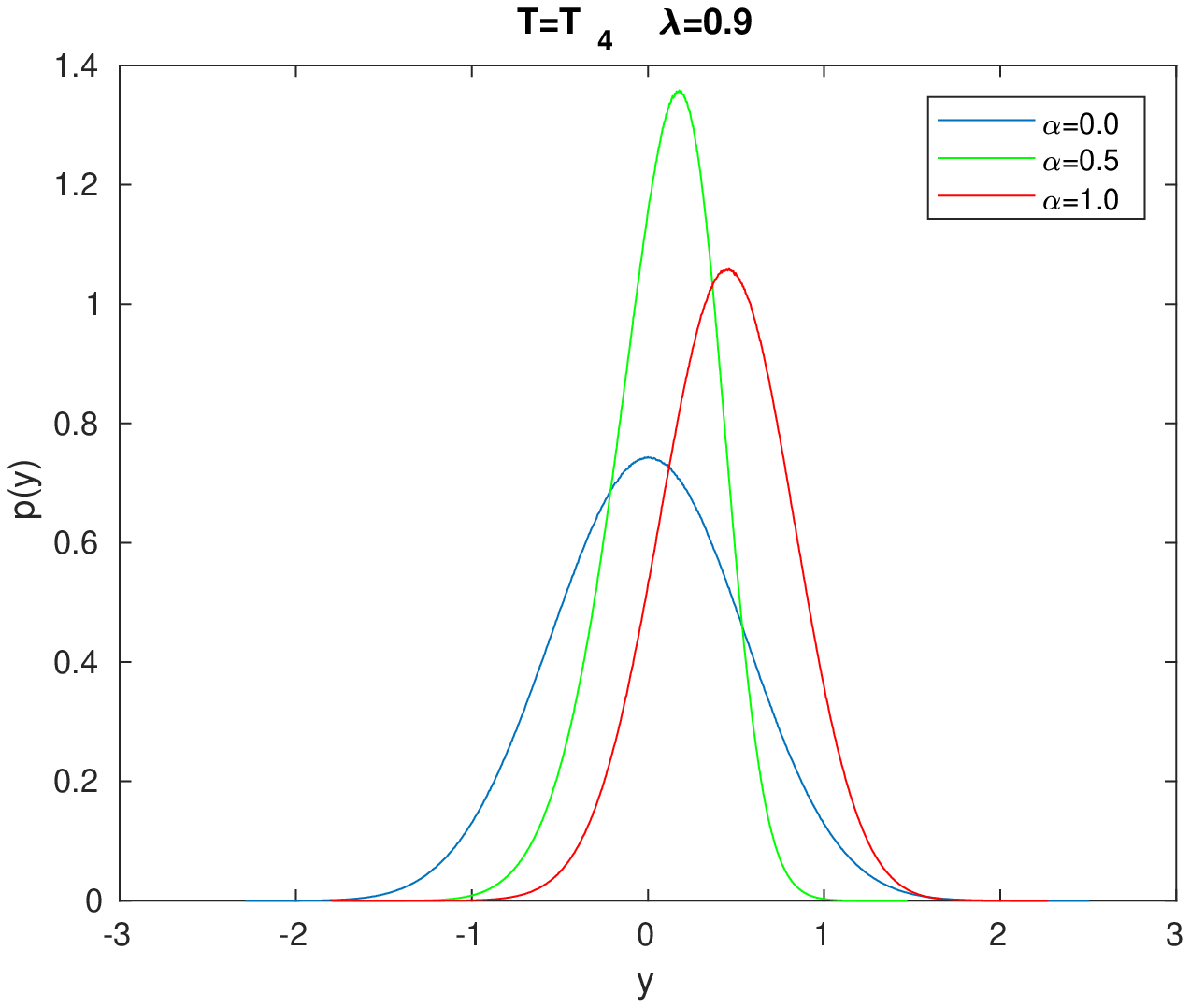}}

        \vspace{4cm}

        \caption{Numerical results of particle velocity distributions
        for driving forces given by diffusively coupled map lattices with couplings $\alpha=0$ (blue), $\alpha=0.5$ (green), and $\alpha=1$ (red), for $T_{2}$ (a),(b); $T_{3}$ (c),(d); $T_{4}$ (e),(f).  The results shown are for $\lambda=e^{-\tau}=0.6$ (a),(c),(e); and $\lambda=e^{-\tau}=0.9$ (b),(d),(f).}
\end{figure}
\begin{figure}[h]
        \subfloat[]{\includegraphics[width=0.4075\textwidth]{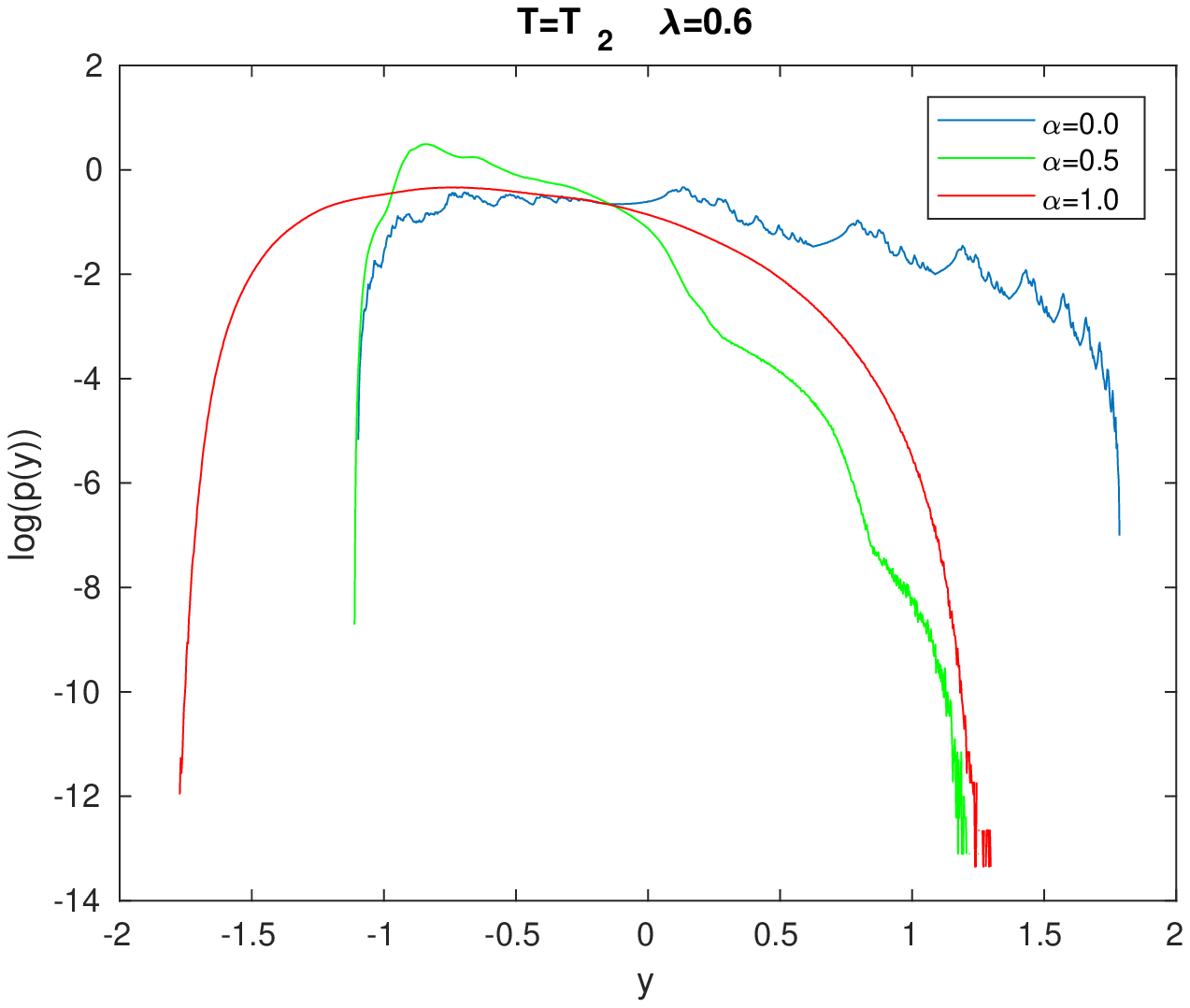}}
        \subfloat[]{\includegraphics[width=0.4075\textwidth]{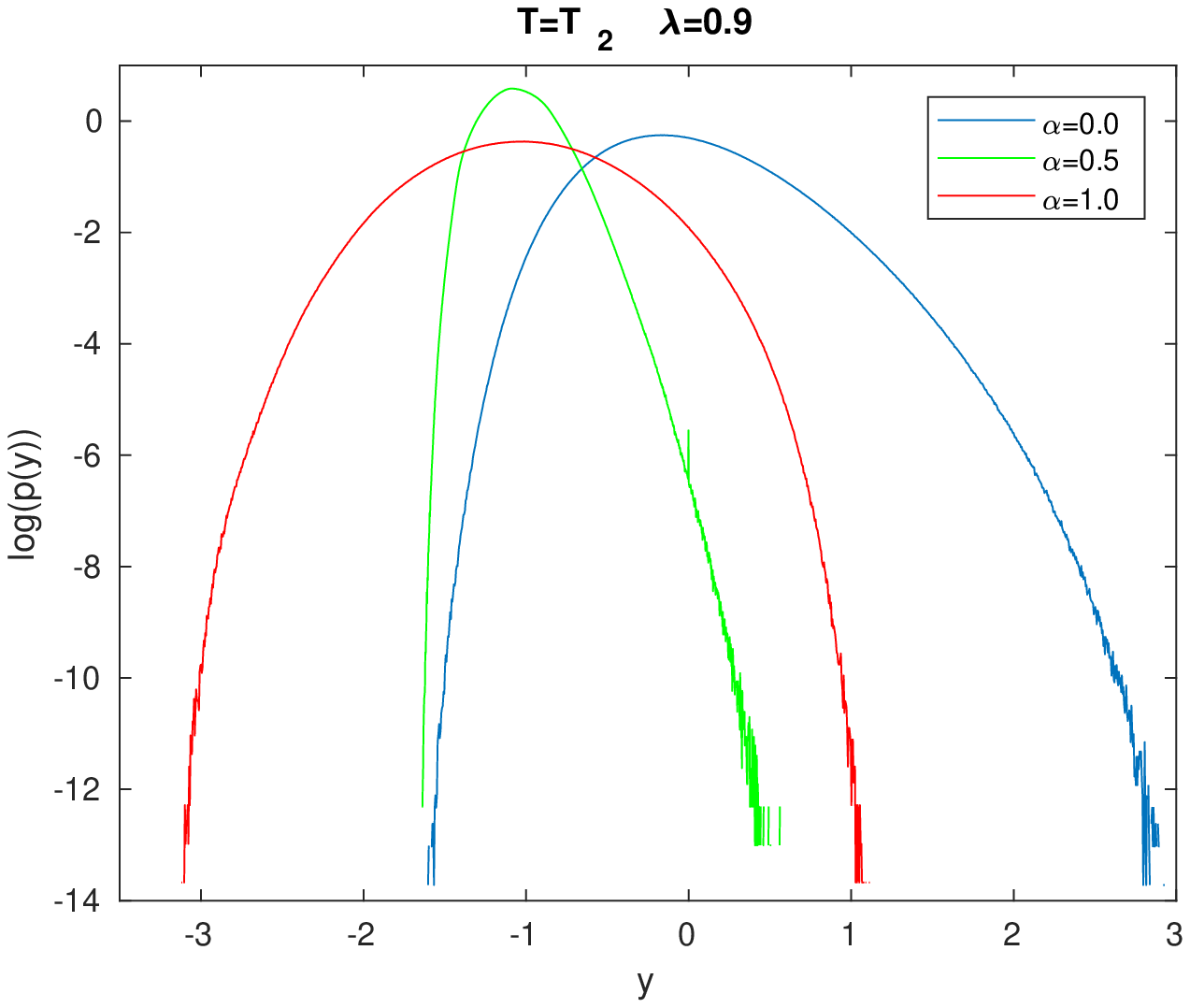}}\\
        \subfloat[]{\includegraphics[width=0.4075\textwidth]{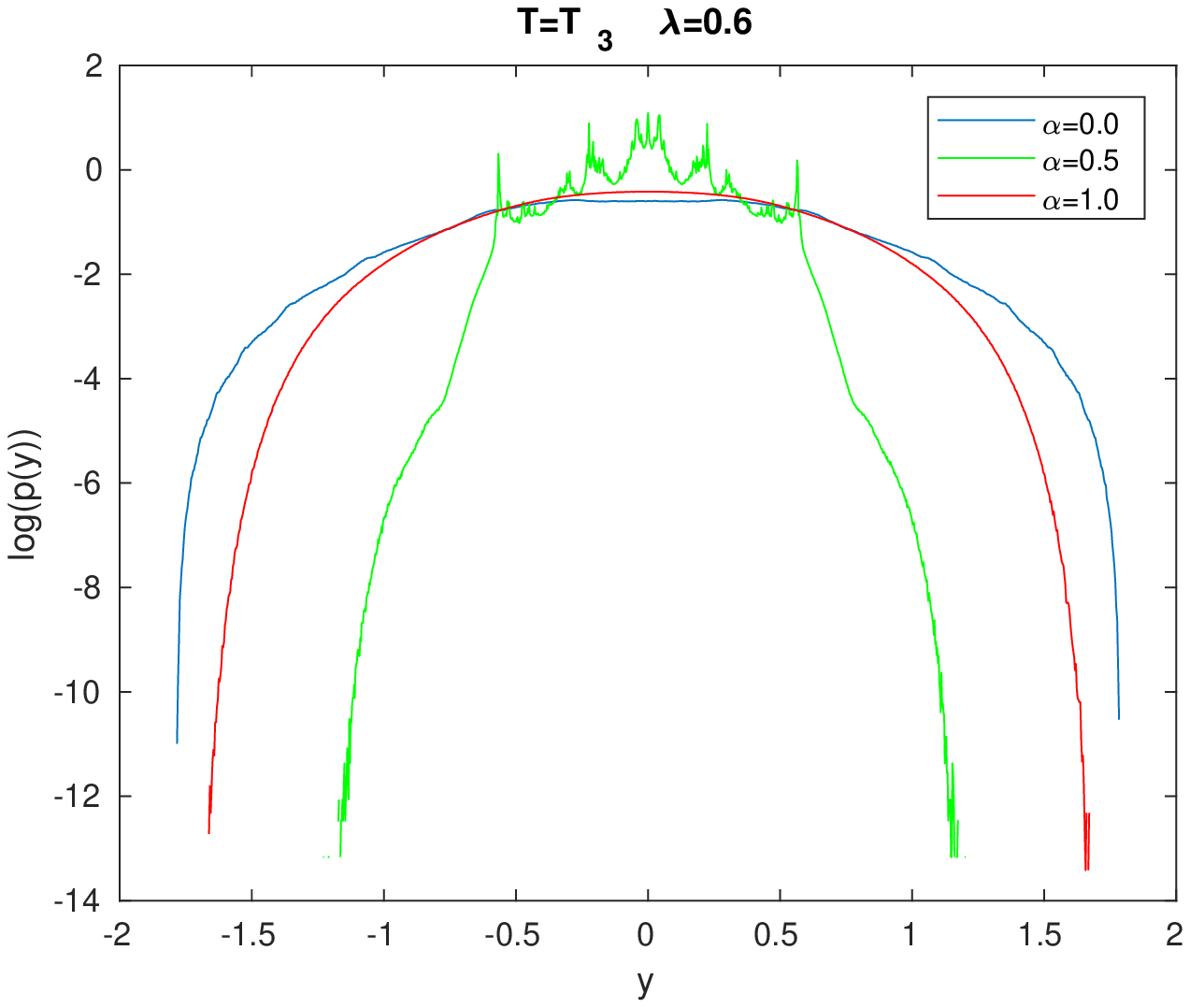}}
        \subfloat[]{\includegraphics[width=0.4075\textwidth]{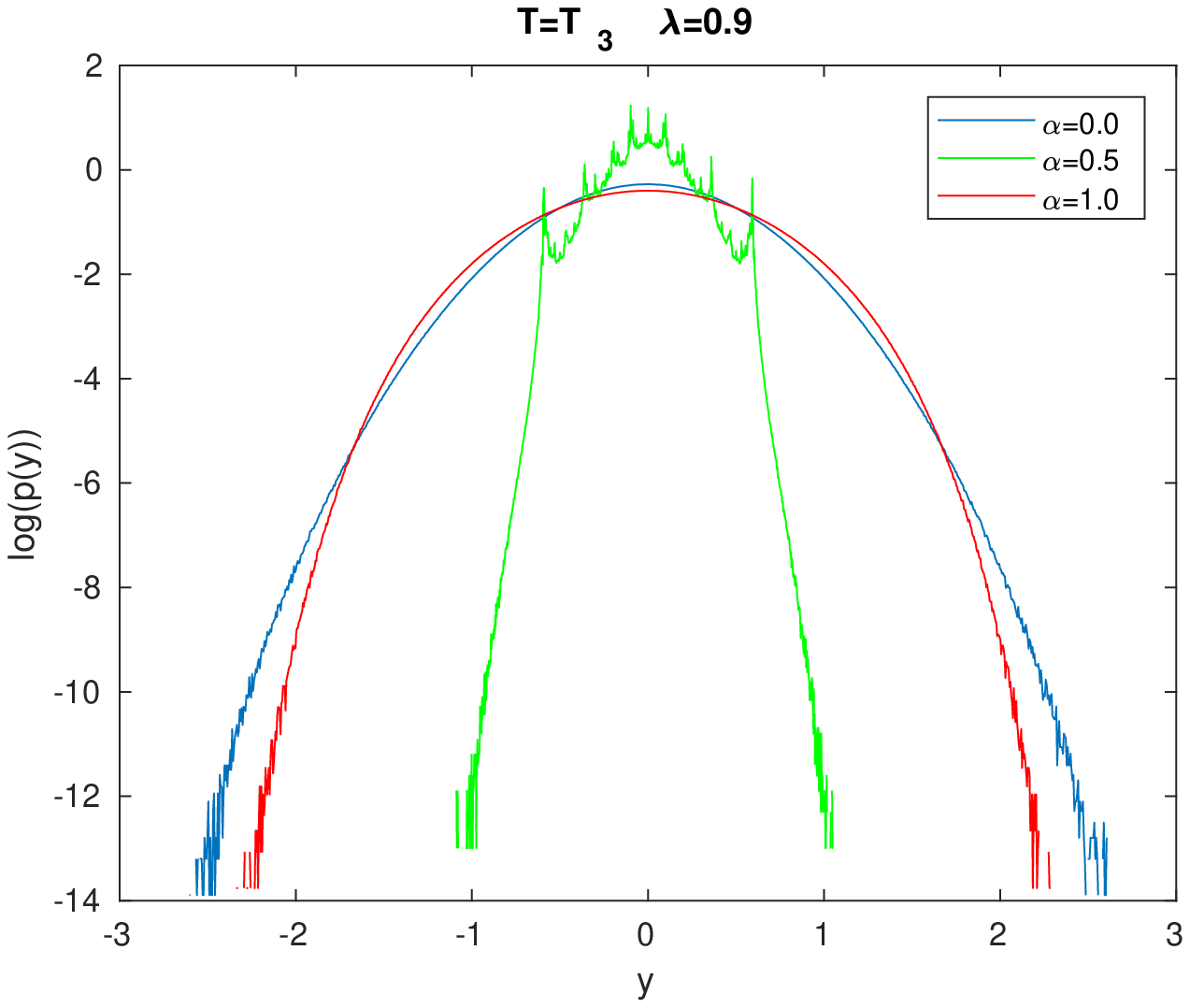}}\\
        \subfloat[]{\includegraphics[width=0.4075\textwidth]{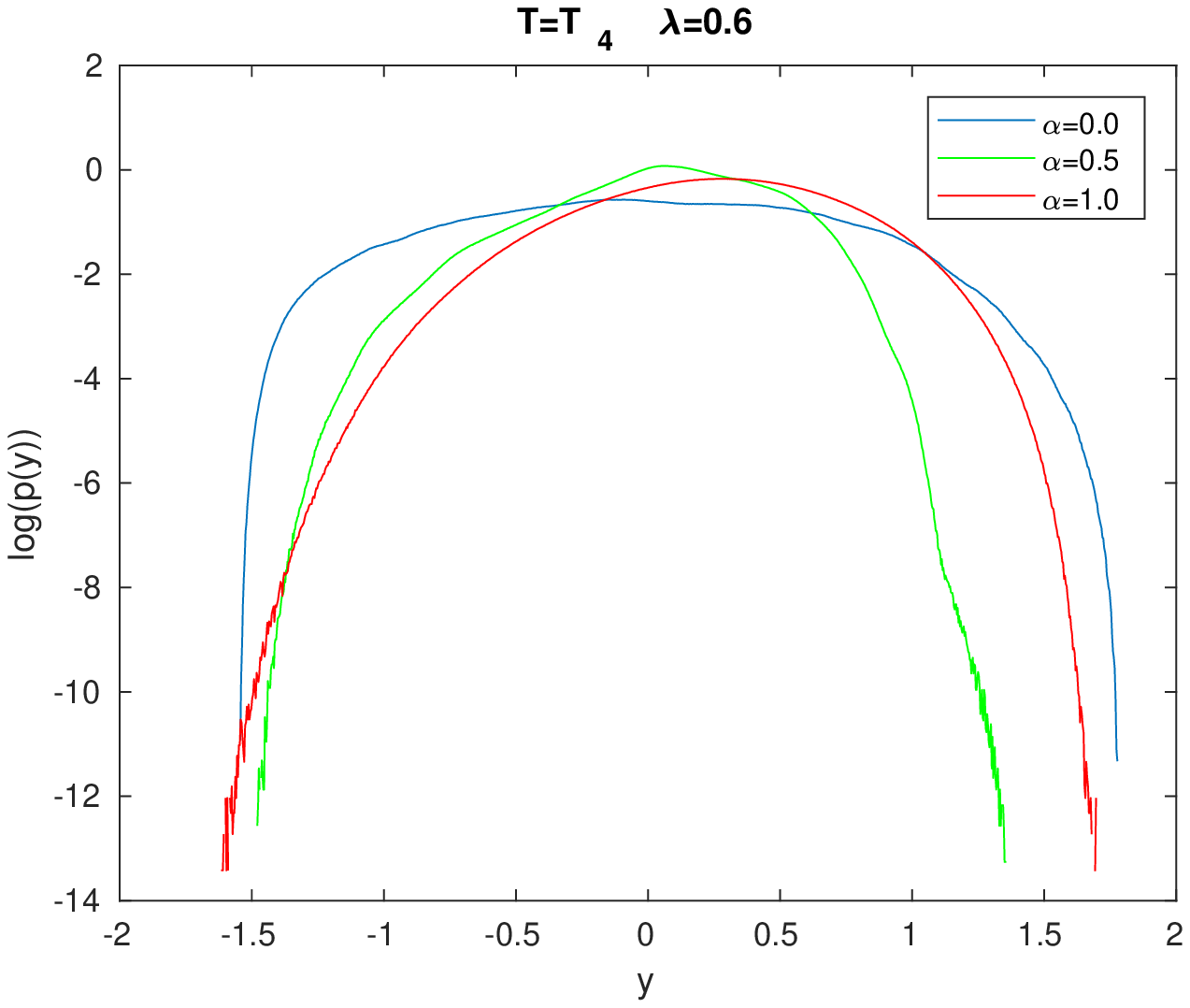}}
        \subfloat[]{\includegraphics[width=0.4075\textwidth]{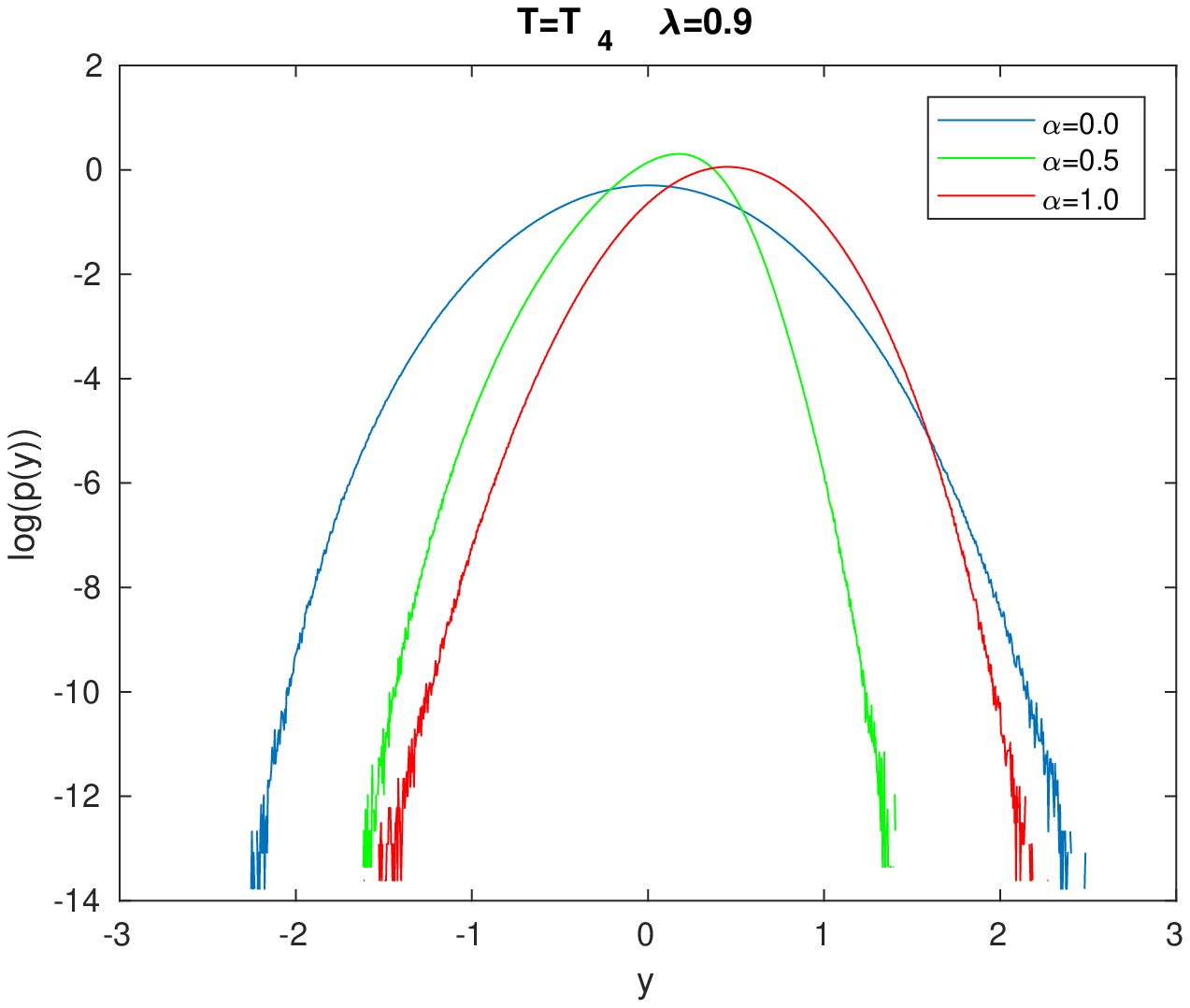}}

        \vspace{4cm}

        \caption{Same as Fig.~4 but in a log-linear plot. In such a plot the Gaussian shape corresponds
        to a parabola.}
\end{figure}
\clearpage

\section{Conclusion}

In this paper we have studied deterministic chaotic analogues of stochastic differential equations, where
the Gaussian white noise of the stochastic differential equation is replaced by deterministic chaotic maps,
evolving on a small time scale $\tau$. We have done a perturbative expansion of the Perron-Frobenius equation
of this problem, and compared terms of equal power in $\sqrt{\tau}$ on both sides of the equation
(a method similar to van Kampen's omega expansion \cite{Kampen} for stochastic systems but here applied in a dynamical systems
context). The result are functional equations which are difficult to solve. However, for Chebyshev maps of
order $N$ we were able to find explicit solutions for arbitrary $N$. The results demonstrate that a
 highly nontrivial mathematical structure arises
where the approach to the Gaussian limit case is different for each of the Chebyshev maps $T_2$, $T_3$, $T_4$.
However, for $N \geq 4$
the leading order correction to the Gaussian is the same for all higher Chebyshev maps $T_N$ with
$N \geq 4$. The extracted polynomial functions in the perturbative expansion of the invariant density can
be regarded as eigenfunctions of a rescaling operator $\tau \to \frac{1}{2} \tau$ that describes via which route the Gaussian limit case is
approached under successive rescaling of the time scale parameter. The leading order
approach to the Gaussian limit case is of order $O(\tau^{1/2})$ for $T_2$ and of order $O(\tau)$ for $N \geq 3$,
though with different polynomial functional forms for $T_3$ and $T_4$. Our results have applications in quantum field theoretical models \cite{Beck7} in which the Gaussian white noise of the Parisis-Wu approach of stochastic quantization is replaced
by a deterministic chaotic dynamics on a small but finite time scale.
{\section{Appendix A}}
Continuing from Eq. (51) in Section IV, {we can evaluate the left-hand side of the inhomogeneous Fokker-Planck equation $(42)$:
\begin{equation}
	\int dx \ xb(x,y,t)=\langle x^{2}\rangle\beta_{1}(y,t)=\frac{1}{2}\left(\frac{2}{\pi}\right)^{1/2}(1-4y^{2})e^{-2y^{2}}
\end{equation}
The stationary solution of Eq. $(42)$ satisfies
\begin{equation}
	\frac{\partial}{\partial y}\alpha(y)+4y\alpha(y)=2\left(\frac{2}{\pi}\right)^{1/2}(1-4y^{2})e^{-2y^{2}}+ const
\end{equation}
In general, all solutions of the first-order linear differential equation
\begin{equation}
	\frac{\partial}{\partial y}\alpha +g(y)\alpha +f(y)=0
\end{equation}
are given by
\begin{equation}
	\alpha(y)=e^{-G(y)}\left[C-\int_{1}f(y)e^{G(y)}\right]
\end{equation}
Here $G$ is an indefinite integral of $g$,
\begin{equation}
	\frac{\partial}{\partial y}G(y)=g(y)
\end{equation}
and $\int_{1}\dots$ denotes an indefinite integral of the argument.  In our case
\begin{equation}
	g(y)=4y
\end{equation}
\begin{equation}
	f(y)=2\left(\frac{2}{\pi}\right)^{1/2}(4y^{2}-1)e^{-2y^{2}} - const
\end{equation}
This yields
\begin{equation}
	G(y)=2y^{2}
\end{equation}
\begin{equation}
	\alpha(y)=e^{-2y^{2}}\left[C-2\left(\frac{2}{\pi}\right)^{1/2}\int_{1}(4y^{2}-1) + const\int_{1}e^{2y^{2}}\right]
\end{equation}
For $C$=const=0 we obtain
\begin{equation}
	\alpha(y)=\left(\frac{2}{\pi}\right)^{1/2}e^{-2y^{2}}\left(2y-\frac{8}{3}y^{3}\right)
\end{equation}
To obtain the second-order correction term, we have to solve the coupled system of equations $(43)$ and $(45)$.  Putting Eq. $(47)$ into Eq. $(45)$ and again using the symmetry of the pre-images of $T$, we can write Eq. $(45)$ as
\begin{equation}
\begin{aligned}
	c(x',y,t) ={} & \sum_{x\in T^{-1}(x')}\frac{1}{\vert T'(x)\vert}\bigg\{c(x,y,t) \\
		& +(1-2x^{2})h(x)\left[\frac{1}{2}\frac{\partial}{\partial y}\beta_{1}(y,t)-\frac{1}{4}\frac{\partial^{2}}{\partial y^{2}}\alpha(y,t)\right]\bigg\}
\end{aligned}
\end{equation}
Similarly as Eq. $(44)$, this equation is solved by the ansatz
\begin{equation}
	c(x,y,t)=h(x)\gamma_{0}(y,t)+xh(x)\gamma_{1}(y,t)
\end{equation}
We obtain $\gamma_{0}(y,t)$ arbitrary and
\begin{equation}
	\gamma_{1}(y,t)=\frac{1}{2}\frac{\partial}{\partial y}\beta_{1}(y,t)-\frac{1}{4}\frac{\partial^{2}}{\partial y^{2}}\alpha(y,t)
\end{equation}
Putting Eqs. $(81)$, $(47)$, and $(50)$ into Eq. $(43)$, we arrive at the following equation for the unknown function $\beta_{0}(y,t)=\beta(y,t)$:
\begin{equation}
\begin{aligned}
	0 ={} & \frac{\partial}{\partial y}(y\beta(y,t))+\frac{1}{4}\frac{\partial^{2}}{\partial y^{2}}\beta(y,t)-\frac{\partial}{\partial t}\beta(y,t) \\
		& +\biggl[\frac{1}{2}y\frac{\partial}{\partial y}+\frac{1}{2}(y^{2}+1)\frac{\partial^{2}}{\partial y^{2}}+\frac{1}{4}y\frac{\partial^{3}}{\partial y^{3}}+\frac{5}{64}\frac{\partial^{4}}{\partial y^{4}} \\
		& -\frac{1}{2}+\frac{\partial}{\partial t}-\frac{1}{2}\frac{\partial^{2}}{\partial t^{2}}\biggr]p_{0}(y,t)+\frac{1}{8}\frac{\partial^{3}}{\partial y^{3}}\alpha(y,t)
\end{aligned}
\end{equation}
In the stationary case, we can use Eqs. $(46)$ and $(79)$ to obtain
\begin{equation}
	0=\frac{\partial}{\partial y}(y\beta(y))+\frac{1}{4}\frac{\partial^{2}}{\partial y^{2}}\beta(y)+\left(\frac{2}{\pi}\right)^{1/2}\left(\frac{64}{3}y^{6}-68y^{4}+46y^{2}-\frac{15}{4}\right)e^{-2y^{2}}
\label{66} \end{equation}
One can easily check that this differential equation is solved by
\begin{equation}
	\beta(y)=\left(\frac{2}{\pi}\right)^{1/2}\left(\frac{32}{9}y^{6}-\frac{31}{3}y^{4}+\frac{15}{2}y^{2}-\frac{37}{48}\right)e^{-2y^{2}}
\end{equation}}
\clearpage
{\section{Appendix B}}
Continuing from Eq. (58) in Section V, and {using Lemma 4 in Appendix C, for $N\geq 3$, $(58)$ reduces to
\begin{equation}
        b(x',y,t)=\sum_{x\in T^{-1}(x')} \frac{1}{|T'(x)|}b(x,y,t)
\end{equation}
which is solved by $b(x,y,t)=h(x)\beta(y,t)$.

The integrated equations reduce for any Chebyshev map with $N\geq 2$ as follows. Since
\begin{equation}
        \int dx \ x \ \varphi(x,y,t)=\langle x \rangle p_{0}(y,t)
\end{equation}
and
\begin{equation}
        \int dx \ x \ a(x,y,t)=\langle x \rangle \alpha(y,t)
\end{equation}
(as per the definition of the first moment, and implied by $(21)$ and $(22)$), we have
\begin{equation} \nonumber
        (26) \Leftrightarrow \langle x \rangle = 0
\end{equation}
\begin{equation} \nonumber
        (27) \Leftrightarrow \frac{\partial}{\partial y}\left(y p_{0}(y,t)\right)+\frac{1}{4}\frac{\partial^{2}}{\partial y^{2}}p_{0}(y,t) - \frac{\partial}{\partial t}p_{0}(y,t) = 0
\end{equation}
(the Fokker-Planck equation), and
\begin{equation} \nonumber
\begin{aligned}
        (28) \Leftrightarrow & \frac{\partial}{\partial y}\int dx \ x \ b(x,y,t) \\
        = & \frac{\partial}{\partial y}\left(y \alpha(y,t)\right)+\frac{1}{4} \frac{\partial^{2}}{\partial y^{2}}\alpha(y,t)-\frac{\partial}{\partial t}\alpha(y,t)
\end{aligned}
\end{equation}
(an inhomogeneous Fokker-Planck equation).  As previously shown, for $N=2$ we have
\begin{equation}
        b(x,y,t)=h(x)\beta_{0}(y,t)+x h(x)\beta_{1}(y,t)
\end{equation}
and $(28)$ becomes $(56)$
\begin{equation} \nonumber
        \frac{\partial}{\partial y}\left(y \alpha(y,t)\right)+\frac{1}{4} \frac{\partial^{2}}{\partial y^{2}}\alpha(y,t)-\frac{\partial}{\partial t}\alpha(y,t) = -\frac{1}{8}\frac{\partial^{3}}{\partial y^{3}}p_{0}(y,t)
\end{equation}
For $N\geq 3$, we have
\begin{equation} \nonumber
	b(x,y,t)=h(x)\beta(y,t)
\end{equation}
hence
\begin{equation} \nonumber
	\int dx \ xb(x,y,t)=\langle x \rangle\beta(y,t)=0
\end{equation}
hence (28) then reads
\begin{equation} \nonumber
	\frac{\partial}{\partial y}(y\alpha(y,t)+\frac{1}{4}\frac{\partial^{2}}{\partial y^{2}}\alpha(y,t)-\frac{\partial}{\partial t}\alpha(y,t)=0
\end{equation}
which is a Fokker-Planck equation.  Now, $(29)$ reduces to
\begin{equation}
\begin{aligned}
        \frac{\partial}{\partial y}\int dx \ xc(x,y,t) = & \frac{\partial}{\partial y}\left(y \beta(y,t)\right)+\frac{1}{4} \frac{\partial^{2}}{\partial y^{2}}\beta(y,t)-\frac{\partial}{\partial t}\beta(y,t) \\
        & +\left(-\frac{1}{2}+\frac{\partial}{\partial t}-\frac{1}{2}\frac{\partial^{2}}{\partial t^{2}}+\frac{1}{2}y\frac{\partial}{\partial y}+\frac{1}{2}(y^{2}+1)\frac{\partial^{2}}{\partial y^{2}}+\frac{1}{4}y\frac{\partial^{3}}{\partial y^{3}}+\frac{1}{64}\frac{\partial^{4}}{\partial y^{4}}\right)e^{-2y^{2}}\sqrt{\frac{2}{\pi}} \\
        = & \frac{\partial}{\partial y}\left(y \beta(y,t)\right)+\frac{1}{4} \frac{\partial^{2}}{\partial y^{2}}\beta(y,t)-\frac{\partial}{\partial t}\beta(y,t)-\sqrt{\frac{2}{\pi}}\left(\frac{7}{4}-10y^{2}+4y^{4}\right)e^{-2y^{2}}
\end{aligned}
\end{equation}
with $\beta=\beta_{0}$ and $\beta_{1}=0$.  $c(x,y,t)$ satisfies Eq. $(18)$, which for $N\geq 3$ can be written as
\begin{equation}
        c(x',y,t) = \sum_{x\in T^{-1}(x')} \frac{1}{\vert T'(x) \vert}\left[c(x,y,t)+\frac{1}{4}h(x)(2x^{2}-1)\frac{\partial^{2}}{\partial y^{2}}\alpha(y,t)-\frac{1}{6}x^{3}h(x)\frac{\partial^{3}}{\partial y^{3}}p_{0}(y,t)\right]
\end{equation}
using Lemmas 2 and 3 in Appendix C and the Fokker-Planck equation for $\alpha(y,t)$.

Using Lemma 4 in Appendix C, $(90)$ even further reduces to
\begin{equation}
        c(x',y,t) = \sum_{x\in T^{-1}(x')} \frac{1}{\vert T'(x) \vert}\left[c(x,y,t)-\frac{1}{6}x^{3}h(x)\frac{\partial^{3}}{\partial y^{3}}p_{0}(y,t)\right]
\end{equation}
Similar to $b(x,y,t)$ for $T_{2}$, we will choose the ansatz for solutions of Eq. $(91)$ as:
\begin{equation}
\begin{aligned}
        c(x,y,t) & =h(x)\gamma_{0}(y,t)+xh(x)\gamma_{1}(y,t) \\
        \Rightarrow & l = h(x')\gamma_{0}(y,t)+x'h(x')\gamma_{1}(y,t)
\end{aligned}
\end{equation}
where $x'=4x^{3}-3x$. This gives
\begin{equation}
\begin{aligned}
        & =\gamma_{0}(y,t)\left(\sum_{x\in T^{-1}(x')}\frac{h(x)}{\vert T'(x) \vert}\right)+\gamma_{1}(y,t)\sum_{x\in T^{-1}(x')}(4x^{3}-3x)\frac{h(x)}{\vert T'(x) \vert} \\
        & = \sum_{x\in T^{-1}(x')}\frac{h(x)}{\vert T'(x) \vert}\left(\gamma_{0}(y,t)+4x^{3}\gamma_{1}(y,t)\right)
\end{aligned}
\end{equation}
by Lemma 3 in Appendix C.
\begin{equation}
\begin{aligned}
        r= & \sum_{x\in T^{-1}(x')} \frac{1}{\vert T'(x) \vert}\left[h(x)\gamma_{0}(y,t)+xh(x)\gamma_{1}(y,t)-\frac{1}{6}x^{3}h(x)\frac{\partial^{3}}{\partial y^{3}}p_{0}(y,t)\right] \\
        = & \sum_{x\in T^{-1}(x')} \frac{h(x)}{\vert T'(x) \vert}\left[\gamma_{0}(y,t)-\frac{1}{6}x^{3}\frac{\partial^{3}}{\partial y^{3}}p_{0}(y,t)\right]
\end{aligned}
\end{equation}
When we let $l=r$ we get $\gamma_{0}(y,t)$ arbitrary and
\begin{equation}
\begin{aligned}
        \gamma_{1}(y,t) = & -\frac{1}{24}\frac{\partial^{3}}{\partial y^{3}}p_{0}(y,t) \\
        = & \left(\frac{2}{\pi}\right)^{1/2}\left(\frac{8}{3}y^{3}-2y\right)e^{-2y^{2}}
\end{aligned}
\end{equation}
Hence
\begin{equation}
        \int dx \ xc(x,y,t)= \langle  x\rangle \gamma_{0}(y,t)+ \langle  x^{2}\rangle \gamma_{1}(y,t)=\left(\frac{2}{\pi}\right)^{1/2}\left(\frac{4}{3}y^{3}-y\right)e^{-2y^{2}}
\end{equation}
(where $\langle  x\rangle = 0$ and $\langle  x^{2}\rangle = \frac{1}{2}$) and
\begin{equation}
\begin{aligned}
        \frac{\partial}{\partial y} \int dx \ xc(x,y,t) = & \left(\frac{2}{\pi}\right)^{1/2}\biggl[\left(4y^{2}-1\right)e^{-2y^{2}}-4y\left(\frac{4}{3}y^{3}-y\right)e^{-2y^{2}}\biggr] \\
        = & \left(\frac{2}{\pi}\right)^{1/2}\left(-\frac{16}{3}y^{4}+8y^{2}-1\right)e^{-2y^{2}}
\end{aligned}
\end{equation}
To this we have to add the $N\geq 4$ contribution (see section VI)
\begin{equation}
        \left(\frac{2}{\pi}\right)^{1/2}\left(4y^{4}-10y^{2}+\frac{7}{4}\right)e^{-2y^{2}}
\end{equation}
Hence, for $N=3$, $\beta(y,t)$ satisfies the inhomogeneous Fokker-Planck equation
\begin{equation}
        \frac{\partial}{\partial y}\left(y\beta(y,t)\right)+\frac{1}{4}\frac{\partial^{2}}{\partial y^{2}}\beta(y,t)- \frac{\partial}{\partial t}\beta(y,t)=\left(\frac{2}{\pi}\right)^{1/2}\left(-\frac{4}{3}y^{4}-2y^{2}+\frac{3}{4}\right)e^{-2y^{2}}
\end{equation}
From the ansatz for the stationary solution
\begin{equation}
        \beta(y)=  \alpha_{0}e^{-2y^{2}}+\alpha_{1}y^{2}e^{-2y^{2}}+\alpha_{2}y^{4}e^{-2y^{2}}
\end{equation}
and using $\int dy \beta(y) =0$ and
\begin{equation} \nonumber
\int e^{-2y^{2}}dy = \left(\frac{2}{\pi}\right)^{1/2}, \ \ \ \ \ \int y^{2}e^{-2y^{2}}dy =\left(\frac{2}{\pi}\right)^{1/2}\frac{1}{4}, \ \ \ \ \ \ \int y^{4}e^{-2y^{2}}dy =\left(\frac{2}{\pi}\right)^{1/2}\frac{3}{16},
\end{equation}
one arrives at
\begin{equation}
\beta(y)=\left(\frac{2}{\pi}\right)^{1/2}\left(\frac{1}{3}y^{4}+\frac{3}{2}y^{2}-\frac{7}{16}\right)e^{-2y^{2}}.
\end{equation}}
\clearpage
\section{Appendix C}
\subsection{Lemma 1}
Let $x'=\cos(\pi N u_{0})=T_{N}(x)$ be the $1^{st}$ iterate of an $N-th$ order Chebyshev polynomial $T_{N}$.  All pre-images of $x'$ are given by
\begin{equation} \nonumber
        x=\cos\left(\pi u_{0}+\frac{2\pi}{N}\cdot
                   j\right)
\end{equation}
where
\begin{equation} \nonumber
\begin{aligned}
        j = & -\frac{N}{2}+1,-\frac{N}{2}+2,\dots , \frac{N}{2} \ \ \ (N \ even) \\
        j = & -\frac{N}{2}+\frac{1}{2},-\frac{N}{2}+\frac{3}{2}, \dots , \frac{N}{2}-\frac{1}{2} \ \ \ (N \ odd)
\end{aligned}
\end{equation}
\textbf{Proof}
\\
Note,
\begin{equation} \nonumber
x' = \cos\left(N\pi u_{0}+2\pi j\right) = \cos\left(N\pi u_{0}\right)
\end{equation}
since $j$ is integer.   For $T_{N}$ there are $N$ pre-images, all characterized by a different $j$. $\square$
\subsection{Lemma 2}
For all $x \in T_{N}^{-1}(x')$ one has
\begin{equation} \nonumber
        \frac{h(x)}{|T_{N}'(x)|}=const_{j}
\end{equation}
independent of the pre-image chosen.
\\
\\
\textbf{Proof}
\\
Note,
\begin{equation} \nonumber
x'= T_{N}(x)=\cos\left(N \arccos (x)\right)=\cos \left(N\pi u_{0}\right)
\end{equation}
if $x=\cos\left(\pi u_{0}\right)$.  This implies
\begin{equation} \nonumber
\begin{aligned}
        |T_{N}'(x)|= & |\sin\left(N \arccos (x) \right)|\cdot N \cdot \frac{1}{\sqrt{1-x^{2}}}=\frac{N\cdot |\sin\left(N\pi u_{0}\right)|}{\sqrt{1-\cos^{2}\left(\pi u_{0}\right)}}  \\
        = & \frac{N\cdot |\sin\left(N\pi u_{0}\right)|}{|\sin\left(\pi u_{0}\right)|}
\end{aligned}
\end{equation}
This further implies
\begin{equation} \nonumber
        h(x)=\frac{1}{\pi\sqrt{1-x^{2}}}=\frac{1}{\pi |\sin\left(\pi u_{0}\right)|} \Rightarrow \frac{h(x)}{|T_{N}'(x)|} = \frac{1}{N\pi |\sin\left(N\pi u_{0}\right)|}
\end{equation}
For each pre-image $u=u_{0}+\frac{2}{N}\cdot j$, one has
\begin{equation} \nonumber
\begin{aligned}
        |\sin\left(N\pi\left(u_{0}+\frac{2}{N}\cdot j\right)\right)| = & |\sin\left(N\pi u_{0} +2\pi j\right)| \\
        = & |\sin\left(N\pi u_{0}\right)|  \ \ \
\end{aligned}
\end{equation}
Hence $\frac{h(x)}{|T_{N}'(x)|}=const_{j}$ (independent of $j$).  In fact, we have for all pre-images $x$
\begin{equation} \nonumber
        \frac{h(x)}{|T_{N}'(x)|}=\frac{1}{N\pi\sqrt{1-\cos^{2}\left(N\pi u_{0}\right)}}=\frac{1}{N\pi\sqrt{1-x'^{2}}}
\end{equation}
which depends on $x'$ only. $\square$
\subsection{Lemma 3}
For all Chebyshev polynomials $T_{N}(x)$, $N \geq 2$, and all $x'$ one has
\begin{equation} \nonumber
        \sum_{x\in T_{N}^{-1}(x')} x=0
\end{equation}
\textbf{Proof}
\\
$N$ even: Trivial (symmetry reason).
\\
$N$ odd: In general,
\begin{equation} \nonumber
\begin{aligned}
        \cos(x)+\cos(y)= & 2\cos\left(\frac{x+y}{2}\right)\cos\left(\frac{x-y}{2}\right) \\
        \Rightarrow & \cos\left(\pi u_{0} +\frac{2\pi}{N}\cdot j\right)+\cos\left(\pi u_{0} -\frac{2\pi}{N}\cdot j\right)=2\cos\left(\pi u_{0}\right)\cos\left(\frac{2\pi}{N}\cdot j\right)
\end{aligned}
\end{equation}
\begin{equation} \nonumber
        \sum_{j=1}^{m}\cos(jx) = \cos\left(\frac{m+1}{2}x\right)\sin\left(\frac{mx}{2}\right)\csc\left(\frac{x}{2}\right) \ \ \ \text{\cite{Grad}}
\end{equation}
where $\csc(x):=\frac{1}{\sin(x)}$.  In our case, $x=\frac{2\pi}{N}$, $m=\frac{N}{2}-\frac{1}{2}$.  This implies
\begin{equation} \nonumber
\begin{aligned}
        \sum_{j=1}^{\frac{N}{2}-\frac{1}{2}} \cos\left(\frac{2\pi}{N}j\right) = & \cos\left(\left(\frac{N}{4}+\frac{1}{4}\right)\frac{2\pi}{N}\right)\sin\left(\left(\frac{N}{4}-\frac{1}{4}\right)\frac{2\pi}{N}\right)\csc\left(\frac{\pi}{N}\right) \\
        = & \cos\left(\frac{\pi}{2}+\frac{\pi}{2N}\right)\sin\left(\frac{\pi}{2}-\frac{\pi}{2N}\right)\csc\left(\frac{\pi}{N}\right) \\
        = & -\sin\left(\frac{\pi}{2N}\right)\cos\left(\frac{\pi}{2N}\right)\csc\left(\frac{\pi}{N}\right) \\
        = & -\frac{1}{2}\sin\left(\frac{\pi}{N}\right)\csc\left(\frac{\pi}{N}\right) = -\frac{1}{2}
\end{aligned}
\end{equation}
Hence
\begin{equation} \nonumber
\begin{aligned}
        \sum_{x \in T_{N}^{-1}(x')} x = & \cos\left(\pi u_{0}\right) + 2\cos\left(\pi u_{0}\right)\sum_{j=1}^{\frac{N}{2}-\frac{1}{2}}\cos\left(\frac{2\pi}{N}j\right) \\
        = & \cos\left(\pi u_{0}\right)\left(1+2\left(\frac{-1}{2}\right)\right) = 0 \ \ \ \square
\end{aligned}
\end{equation}
\subsection{Lemma 4}
For all Chebyshev polynomials $T_{N}(x)$, $N \geq 3$, and all $x'$ one has
\begin{equation} \nonumber
        \sum_{x \in T_{N}^{-1}(x')} T_{2}(x) =0
\end{equation}
\textbf{Proof}
\\
$x = \cos\left(\pi u\right) \Rightarrow T_{2}(x) = \cos\left(2\pi u\right)$
\\
Instead of $\cos\left(\pi u_{0}+\frac{2\pi}{N}j\right)$ we thus have to consider $\cos\left(2\pi u_{0}+\frac{4\pi}{N}j\right)$ in the previous proof.
\\
N odd:
\begin{equation} \nonumber
\begin{aligned}
        & \cos\left(2\pi u_{0}+\frac{4\pi}{N}j\right)+\cos\left(2\pi u_{0}-\frac{4\pi}{N}j\right) \\
        & = 2\cos\left(2\pi u_{0}\right)\cos\left(\frac{4\pi}{N}j\right)
\end{aligned}
\end{equation}
The same proof as for Lemma 3 applies here, replacing $\frac{2\pi}{N}$ with $\frac{4\pi}{N}$.
\begin{equation} \nonumber
        \sum_{x \in T_{N}^{-1}(x')} T_{2}(x) = \cos\left(2\pi u_{0}\right)+2\cos\left(2\pi u_{0}\right)\left(\frac{-1}{2}\right)=0
\end{equation}
N even:
\begin{equation} \nonumber
\begin{aligned}
        \sum_{x \in T_{N}^{-1}(x')} T_{2}(x) = & \\
        = & \cos\left(2\pi u_{0}\right) + \cos\left(2\pi u_{0}+2\pi\right)+\sum_{j=1}^{\frac{N}{2}-1}\left(\cos\left(2\pi u_{0}+\frac{4\pi}{N}j\right)+\cos\left(2\pi u_{0}-\frac{4\pi}{N}j\right)\right)
\end{aligned}
\end{equation}
where $\cos\left(2\pi u_{0}\right)$ corresponds to $j=0$ and $\cos\left(2\pi u_{0}+2\pi\right)$ corresponds to $j=\frac{N}{2}$.  This continues as
\begin{equation} \nonumber
\begin{aligned}
        = & 2\cos\left(2\pi u_{0}\right)+2\sum_{j=1}^{\frac{N}{2}-1} \left(\cos\left(2\pi u_{0}\right)\cos\left(\frac{4\pi}{N}j\right)\right) \\
        = & 2\cos\left(2\pi u_{0}\right)\left(1+\sum_{j=1}^{\frac{N}{2}-1}\cos\left(\frac{4\pi}{N}j\right)\right)
\end{aligned}
\end{equation}
$m=\frac{N}{2}-1$, $x=\frac{4\pi}{N}$ $\Rightarrow$
\begin{equation} \nonumber
\begin{aligned}
        \cos\left(\frac{m+1}{2}\right)\sin\left(\frac{mx}{2}\right)\csc\left(\frac{x}{2}\right) = & \cos\left(\frac{N}{4}\frac{4\pi}{N}\right)\sin\left(\left(\frac{N}{4}-\frac{1}{2}\right)\frac{4\pi}{N}\right)\frac{1}{\sin\left(\frac{2\pi}{N}\right)} \\
        = & -1 \cdot \sin\left(\pi - \frac{2\pi}{N}\right)\frac{1}{\sin\left(\frac{2\pi}{N}\right)} \\
        = & -1 \\
        \Rightarrow & \sum_{x \in T_{N}^{-1}(x')} T_{2}(x) = 2\cos\left(2\pi u_{0}\right)(1-1) = 0 \ \ \ \square
\end{aligned}
\end{equation}
\subsection{Lemma 5}

(Generalisation) Let $T_{N}$, $T_{M}$ be suitable Chebyshev polynomials with $N>M \geq 1$.  One has
\begin{equation} \nonumber
        \sum_{x \in T_{N}^{-1}(x')} T_{M}(x) = 0
\end{equation}
\textbf{Proof}
\\
$x=\cos(\pi u) \Rightarrow T_{M}(x) = \cos(M\pi u)$.  Have to consider $\cos\left(M\pi u_{0} +\frac{2\pi M}{N}j\right)$.  The proof is similar as before. $\square$

\end{document}